\documentclass[12pt]{article}
\usepackage{amsfonts}
\usepackage{sw20elba}
\usepackage[T1]{fontenc}

\newtheorem{theorem}{Theorem}

\newtheorem{corollary}[theorem]{Corollary}

\newtheorem{lemma}[theorem]{Lemma}

\newtheorem{proposition}[theorem]{Proposition}

\newenvironment{proof}[1][Proof]{\noindent\textbf{#1.} }{\ \rule{0.5em}{0.5em}}
\input{tcilatex}
\begin{document}

\title{Geometric Quantization of Algebraic Reduction\thanks{%
Revised version of reference \cite{sniatycki06}}}
\author{J\k{e}drzej \'{S}niatycki \\
Department of Mathematics and Statistics\\
University of Calgary\\
Calgary, Alberta, Canada}
\date{}
\maketitle

\begin{abstract}
\noindent We develop an analogue of geometric quantization for Poisson
algebras obtained by algebraic reduction and apply it to generalize the
results of Guillemin and Sternberg on commutation of quantization and
reduction.

\noindent
\_\_\_\_\_\_\_\_\_\_\_\_\_\_\_\_\_\_\_\_\_\_\_\_\_\_\_\_\_\_\_\_\_\_\_\_\_\_

\end{abstract}

\section{Introduction}

For a proper and free Hamiltonian action of a Lie group $G$ on a symplectic
manifold $(P,\omega )$ the space $P/G$ of $G$-orbits in $P$ is a manifold.
If the action of $G$ on $P$ is proper but not free, $P/G$ is a stratified
space. Its geometry is completely determined by the ring $C^{\infty }(P)^{G}$
of $G$-invariant smooth functions on $P$. The description of the orbit space
in terms $C^{\infty }(P)^{G}$, called \emph{singular reduction}, was
initiated in examples by Cushman \cite{cushman83} and abstractly formulated
by Arms, Cushman and Gotay \cite{arms-cushman-gotay}. If the action of $G$
is not proper, then smooth $G$-invariant functions need not separate orbits
and singular reduction is inadequate.

In 1983 Weinstein and I introduced \emph{algebraic reduction} of the zero
level of the momentum map which gives an interesting Poisson algebra even
for an improper action \cite{sniatycki-weinstein}. The extension of
algebraic reduction to non-zero levels of the momentum map was given by
Kimura \cite{kimura93}, Wilbour \cite{wilbour93} and Arms \cite{arms}. For a
free and proper action, all reduction techniques are equivalent. However, if
the action is not free, the Poisson algebra of algebraic reduction need not
be an algebra of functions even if the group $G$ is compact \cite%
{arms-gotay-jennings}.

In 1982 Guillemin and Sternberg proved that quantization and reduction
commute if the action of $G$ is free and proper, $(P,\omega )$ is a compact K%
\"{a}hler manifold \cite{guillemin-sternberg82}. Generalizations of this
result to non-free actions have been studied by many authors, see \cite%
{meinrenken-sjamaar} and references quoted there.

Partially motivated by algebraic reduction, Huebschmann developed a theory
of quantization of abstract Poisson algebras \cite{huebschmann90}. Because
of its generality, Huebschmann's theory is rather difficult to apply to
special cases. In particular, for Poisson algebras of algebraic reduction,
it does not assume \`{a} priori relations between the quantization structure
of the original symplectic manifold and the quantization structure of the
reduced Poisson algebra. Hence, it does not facilitate a study of
commutativity of quantization and algebraic reduction.

The aim of this paper is to develop a theory of quantization of Poisson
algebras obtained by algebraic reduction in terms of a quantization
structure of the original symplectic manifold. Our results provide a
justification for an ad hoc quantization of the example discussed in \cite%
{sniatycki-weinstein}.

We discuss two special cases. In the first case, we consider quantization in
terms of a real polarization and assume that the momentum map is constant
along the leaves of the polarization. We show that there is a natural
duality between the space of generalized invariant sections of the original
quantization and the representation space of the reduced quantization. We
illustrate this case with an example of the action of $G=\mathbb{R}$ on the
phase space of a particle with one degree of freedom generated by its
kinetic energy.

In the second case, we consider quantization of a group action on a K\"{a}%
hler manifold $P$ (not necessarily compact) leading to a unitary
quantization representation $\mathcal{U}$ of $G$ on a Hilbert space $%
\mathcal{H}$. For a quantizable co-adjoint orbit $O$, leading to an
irreducible unitary representation $\mathcal{U}_{O}$ of $G$, we use
quantization of algebraic reduction at $O$ to construct a projection
operator onto the closed subspace of $\mathcal{H}$ on which the
representation $\mathcal{U}$ is equivalent to a multiple (possibly infinite)
of $\mathcal{U}_{O}$. This is a generalization of the results on
commutativity of quantization and reduction obtained in a compact case by
Guillemin and Sternberg \cite{guillemin-sternberg82} and Meinrenken and
Sjamaar \cite{meinrenken-sjamaar}.

It should be noted that proofs given here are quite simple. Yet, the
obtained results incorporate some of the results obtained in the program of
research on commutativity of quantization and reduction and generalize them
to non-compact setting. This suggests that algebraic reduction is a useful
tool in studying the decomposition of a quantization representation into
irreducible components.

\section{Symplectic manifolds and algebraic reduction}

In this section, we review elements of symplectic geometry and algebraic
reduction. The sign convention used here follows reference \cite{sniatycki80}%
.

\subsection{Symplectic manifolds}

Let $(P,\omega )$ be a symplectic manifold. For every $f\in C^{\infty }(P)$,
we denote by $X_{f}$ the Hamiltonian vector field of $f$ defined by 
\begin{equation}
X_{f}%
{\mbox{$ \rule {5pt} {.5pt}\rule {.5pt} {6pt} \, $}}%
\omega =-df,  \label{1}
\end{equation}%
where $%
{\mbox{$ \rule {5pt} {.5pt}\rule {.5pt} {6pt} \, $}}%
$ denotes the left interior product of vectors and forms. Non-degeneracy of $%
\omega $ ensures that $X_{f}$ is well defined. Hamiltonian vector fields
preserve $\omega $ because 
\[
\pounds _{X_{f}}\omega =X_{f}%
{\mbox{$ \rule {5pt} {.5pt}\rule {.5pt} {6pt} \, $}}%
d\omega +d(X_{f}%
{\mbox{$ \rule {5pt} {.5pt}\rule {.5pt} {6pt} \, $}}%
\omega )=0\text{,} 
\]%
where $\pounds $ denotes Lie derivative. The commutator of two Hamiltonian
vector fields $X_{f_{1}}$ and $X_{f_{2}}$ satisfies a relation 
\begin{equation}
\lbrack X_{f_{1}},X_{f_{2}}]%
{\mbox{$ \rule {5pt} {.5pt}\rule {.5pt} {6pt} \, $}}%
\omega =-X_{f_{1}}f_{2}=-d(\omega (X_{f_{1}},X_{f_{2}})).  \label{1.1}
\end{equation}

The symplectic form $\omega $ on $P$ gives rise to a Poisson bracket on $%
C^{\infty }(P)$ defined by 
\begin{equation}
\{f_{1},f_{2}\}=-\omega (X_{f_{1}},X_{f_{2}}).  \label{2}
\end{equation}%
Equation (\ref{1.1}) implies that 
\[
\lbrack X_{f_{1}},X_{f_{2}}]%
{\mbox{$ \rule {5pt} {.5pt}\rule {.5pt} {6pt} \, $}}%
\omega =-X_{\{f_{1},f_{2}\}},
\]%
and 
\begin{equation}
\{f_{1},f_{2}\}=-X_{f_{1}}f_{2}.  \label{3}
\end{equation}%
The Poisson bracket satisfies Jacobi's identity%
\[
\{f_{1},\{f_{2},f_{3}\}\}+\{f_{2},\{f_{3},f_{1}\}\}+\{f_{3},\{f_{1},f_{2}\}%
\}=0,
\]%
and Leibniz' rule 
\[
\{f_{1},f_{2}f_{3}\}=f_{2}\{f_{1},f_{3}\}+\{f_{1},f_{2}\}f_{3}
\]%
for all $f_{1},f_{2},f_{3}\in C^{\infty }(P)$. The ring $C^{\infty }(P)$
endowed with the Poisson bracket given by equation (\ref{2}) is called the
Poisson algebra of $(P,\omega )$. The map $f\mapsto X_{f}$ is an
anti-homomorphism of the Lie algebra structure of the Poisson algebra $%
C^{\infty }(P)$ to the Lie algebra of smooth vector fields on $P$.

\subsection{Momentum map}

Let $G$ be a connected Lie group with Lie algebra $\mathfrak{g}$ and 
\begin{equation}
\Phi :G\times P\rightarrow P:(g,p)\mapsto \Phi _{g}(p)\equiv gp  \label{4}
\end{equation}%
be an action of $G$ on $P$. The action $\Phi $ is free if $gp=p$ implies
that $g$ is the identity in $G$. The action is proper if, for every
convergent sequence $(p_{n})$ in $P$ and a sequence $(g_{n})$ in $G$, such
that the sequence $(g_{n}p_{n})$ converges, the sequence $(g_{n})$ has a
convergent subsequence $(g_{n_{k}})$ and $\lim_{n\rightarrow \infty
}(g_{n}p_{n})=(\lim_{k\rightarrow \infty }g_{n_{k}})(\lim_{n\rightarrow
\infty }p_{n}).$

For each $\xi \in \mathfrak{g}$, the action of $\exp t\xi $ on $P$ gives
rise to a vector field $X^{\xi }$ such that%
\[
(X^{\xi }f)(x)=\frac{d}{dt}f((\exp t\xi )x)_{\mid t=0} 
\]%
for every $f\in C^{\infty }(P)$. The map $\xi \mapsto X^{\xi }$ is an
anti-homomorphism of $\mathfrak{g}$ to the Lie algebra of vector fields on $%
P $.

A momentum map for the action $\Phi $ of $G$ on $P$ is a map $J:P\rightarrow 
\mathfrak{g}^{\ast }$ such that $X^{\xi }$ is the Hamiltonian vector field
of $J_{\xi }=\langle J\mid \xi \rangle $ for every $\xi \in \mathfrak{g}$,
where $\langle J\mid \xi \rangle $ is the evaluation of $J$ on $\xi $. In
other words, 
\[
X^{\xi }=X_{J_{\xi }} 
\]%
for all $\xi \in \mathfrak{g}$. An action of $G$ on $P$ is \emph{Hamiltonian}
if it admits an $Ad^{\ast }$-equivariant momentum map $J:P\rightarrow 
\mathfrak{g}^{\ast }$. The condition that $J$ is $Ad^{\ast }$-equivariant
means that 
\begin{equation}
\Phi _{g}^{\ast }J=Ad_{g}^{\ast }J  \label{5}
\end{equation}%
for every $g\in G$. For $\xi ,\zeta \in \mathfrak{g}$, equation (\ref{5})
implies that 
\[
\langle J\circ \Phi _{\exp t\xi }\mid \zeta \rangle =\langle Ad_{\exp t\xi
}^{\ast }J\mid \zeta \rangle =\langle J\mid Ad_{\exp (-t\xi )}\zeta \rangle
. 
\]%
Differentiating with respect to $t$ and setting $t=0$ we get%
\begin{eqnarray*}
\{J_{\xi },J_{\zeta }\} &=&-X_{J_{\xi }}J_{\zeta }=-\frac{d}{dt}\langle
J\circ \Phi _{\exp t\xi }\mid \zeta \rangle _{\mid t=0}=-\frac{d}{dt}\langle
J\mid Ad_{\exp (-t\xi )}\zeta \rangle _{\mid t=0} \\
&=&-J_{-[\xi ,\zeta ]}=J_{[\xi ,\zeta ]}.
\end{eqnarray*}%
Thus, for a Hamiltonian action of $G$ on $P$, the map $\mathfrak{g}%
\rightarrow C^{\infty }(P):\xi \mapsto J_{\xi }$ is a Lie algebra
homomorphism.

\subsection{Co-adjoint orbits}

Let $O\subset \mathfrak{g}^{\ast }$ be a co-adjoint orbit of $G$. If $\mu
\in O$, then 
\[
O=\{Ad_{g}^{\ast }\mu \mid g\in G\}.
\]%
For each $\xi \in \mathfrak{g}$, we denote by $X_{O}^{\xi }$ the vector
field on $O$ induced by the co-adjoint action of $\exp t\xi $ on $\mathfrak{g%
}^{\ast }$. The symplectic form $\omega _{O}$ of $O$ is given by 
\[
\omega _{O}(X_{O}^{\xi }(\mu ),X_{O}^{\zeta }(\mu ))=-\langle \mu \mid
\lbrack \xi ,\zeta ]\rangle 
\]%
for all $\xi ,\zeta \in \mathfrak{g}$ and $\mu \in O$. Let $%
J_{O}:O\rightarrow \mathfrak{g}^{\ast }$ denote the inclusion map, and $%
J_{O\zeta }\in C^{\infty }(O)$ be given by $J_{O\zeta }(\mu )=\langle
J_{O}(\mu )\mid \zeta \rangle =\langle \mu \mid \zeta \rangle $ for each $%
\mu \in O$ and $\zeta \in \mathfrak{g}$. Then, 
\begin{eqnarray*}
(X_{O}^{\xi }J_{O\zeta })(\mu ) &=&X_{O}^{\xi }(\langle J_{O}\mid \zeta
\rangle )(\mu )=\frac{d}{dt}\langle J_{O}(Ad_{\exp t\xi }^{\ast }(\mu ))\mid
\zeta \rangle _{\mid t=0} \\
&=&\frac{d}{dt}\langle Ad_{\exp t\xi }^{\ast }\mu \mid \zeta \rangle _{\mid
t=0}=\frac{d}{dt}\langle \mu \mid Ad_{\exp (-t\xi )}\zeta \rangle _{\mid t=0}
\\
&=&\langle \mu \mid -[\xi ,\zeta ]\rangle =\omega _{O}(X_{O}^{\xi }(\mu
),X_{O}^{\zeta }(\mu )).
\end{eqnarray*}%
Hence, 
\[
X_{O}^{\zeta }%
{\mbox{$ \rule {5pt} {.5pt}\rule {.5pt} {6pt} \, $}}%
\omega =-dJ_{O\zeta },
\]%
which shows that $J_{O}$ is a momentum map for the co-adjoint action of $G$
on $O$.

\bigskip

\subsection{Algebraic reduction at $0\in \mathfrak{g}^{\ast }$}

Let $(P,\omega )$ be a symplectic manifold, $\Phi :G\times P\rightarrow P$ a
Hamiltonian action on $P$ of a connected group $G$, and $J:P\rightarrow 
\mathfrak{g}^{\ast }$ an $Ad^{\ast }$-equivariant momentum map for this
action. Let $\mathcal{J}$ be the associative ideal in $C^{\infty }(P)$
generated by components of $J.$ In other words, if $(\xi _{1},...,\xi _{k})$
is a basis of the Lie algebra $\mathfrak{g}$ of $G$, then 
\begin{equation}
f\in \mathcal{J}~~\Longleftrightarrow ~~f=\sum_{j=1}^{k}f_{j}J_{\xi _{j}}
\label{6}
\end{equation}%
for some functions $f_{1},...,f_{k}\in C^{\infty }(P)$. Let $C^{\infty }(P)/%
\mathcal{J}$ be the space of equivalence classes of smooth functions on $P$
modulo $\mathcal{J}$. For every $f\in C^{\infty }(P)$, the class of $f$ in $%
C^{\infty }(P)/\mathcal{J}$ is 
\[
\lbrack f]=\left\{ f+\sum_{j=1}^{k}f_{j}J_{\xi _{j}}\mid f_{1},...,f_{k}\in
C^{\infty }(P)\right\} . 
\]

Since the momentum mapping $J$ is $Ad^{\ast }$-equivariant, it follows that 
\[
\Phi _{g}^{\ast }\left( \sum_{j=1}^{k}f_{j}J_{\xi _{j}}\right)
=\sum_{j=1}^{k}(\Phi _{g}^{\ast }f_{j})J_{Ad_{g}(\xi _{j})} 
\]%
for every $g\in G$. Hence, $\mathcal{J}$ is $G$-invariant, and there is an
induced action $\tilde{\Phi}^{\ast }$ of $G$ on $C^{\infty }(P)/\mathcal{J}$
such that 
\[
\tilde{\Phi}_{g}^{\ast }[f]=[\Phi _{g}^{\ast }f] 
\]%
for every $f\in C^{\infty }(P)$ and $g\in G$. We denote by $(C^{\infty }(P)/%
\mathcal{J})^{G}$ the space of $G$-invariant elements of $C^{\infty }(P)/%
\mathcal{J}$. If $[f]$ is $G$-invariant, $[X_{J_{\xi }}f]=[f]$ for all $\xi
\in \mathfrak{g}$. By assumption, $G$ is connected so that 
\begin{equation}
(C^{\infty }(P)/\mathcal{J})^{G}=\{[f]\in C^{\infty }(P)/\mathcal{J}\mid
X_{J_{\xi }}f\in \mathcal{J}\text{ for all }\xi \in \mathfrak{g}\}.
\label{8}
\end{equation}%
For $[f_{1}],[f_{2}]\in C^{\infty }(P)/\mathcal{J}$, the product $%
[f_{1}][f_{2}]$ is defined by 
\[
\lbrack f_{1}][f_{2}]=[f_{1}f_{2}]. 
\]%
If $[f_{1}]$ and $[f_{2}]$ are $G$-invariant, then $X_{J_{\xi
}}(f_{1}f_{2})=(X_{J_{\xi }}f_{1})f_{2}+f_{1}(X_{\xi }f_{2})\in \mathcal{J}$
for every $\xi \in \mathfrak{g}$. Hence, $(C^{\infty }(P)/\mathcal{J})^{G}$
is an associative subalgebra of $C^{\infty }(P)/\mathcal{J}$. Moreover, $%
(C^{\infty }(P)/\mathcal{J})^{G}$ is also a Poisson algebra with a Poisson
bracket defined by 
\begin{equation}
\{[f_{1}],[f_{2}]\}=[\{f_{1},f_{2}\}].  \label{9}
\end{equation}%
In order to see that the Poisson bracket on $(C^{\infty }(P)/\mathcal{J}%
)^{G} $ is well defined observe that, if $f_{1}=fJ_{\xi }\in \mathcal{J}$
and $[f_{2}]\in (C^{\infty }(P)/\mathcal{J})^{G}$, then 
\[
\lbrack \{f_{1},f_{2}\}]=[\{fJ_{\xi },f_{2}\}]=[f\{J_{\xi
},f_{2}\}+\{f,f_{2}\}J_{\xi }]=0 
\]%
because $\{J_{\xi },f_{2}\}=-X_{J_{\xi }}f_{2}\in \mathcal{J}$ by equation (%
\ref{8}). Moreover, 
\[
X_{J_{\xi }}\{f_{1},f_{2}\}=\{X_{J_{\xi }}f_{1},f_{2}\}+\{f_{1},X_{J_{\xi
}}f_{2}\} 
\]%
for every $f_{1},f_{2}\in C^{\infty }(P)$ and $\xi \in \mathfrak{g}.$ If $%
[f_{1}],[f_{2}]\in (C^{\infty }(P)/\mathcal{J})^{G}$ then $X_{J_{\xi }}f_{1}$
and $X_{J_{\xi }}f_{2}$ are in $\mathcal{J}$ and the above argument implies
that $X_{\xi }\{f_{1},f_{2}\}\in \mathcal{J}$. Hence, equation (\ref{9})
gives a well defined bracket on $(C^{\infty }(P)/\mathcal{J})^{G}$.
Linearity, antisymmetry, Leibniz' rule and Jacobi's identity of this bracket
follow from the corresponding properties of the Poisson bracket on $%
C^{\infty }(P)$.

The Poisson algebra $(C^{\infty }(P)/\mathcal{J})^{G}$ is the Poisson
algebra of algebraic reduction of $(P,\omega )$ at the zero value of the
momentum map $J$. For a free and proper Hamiltonian action of $G$ on $%
(P,\omega )$, the orbits space $J^{-1}(0)/G$ is a symplectic manifold,
called a Marsden-Weinstein reduced space, and the the Poisson algebra $%
(C^{\infty }(P)/J)^{G}$ is isomorphic to the Poisson algebra $C^{\infty
}(J^{-1}(0)/G)$, \cite{sniatycki-weinstein}.

\subsection{Algebraic reduction at non-zero orbits.}

Consider now a non-zero co-adjoint orbit $O=\{{Ad}_{g}^{\ast }\mu \mid g\in
G\}$. We denote by $G_{\mu }$ the isotropy group of $\mu $. In other words, 
\[
G_{\mu }=\{g\in G\mid Ad_{g}^{\ast }\mu =\mu \}. 
\]%
If the action of $G$ on $P$ is free and proper, then $J^{-1}(\mu )/G_{\mu }$
is a symplectic manifold, \cite{marsden-weinstein}. Moreover, there is a
bijection between $J^{-1}(O)/G$ and $J^{-1}(\mu )/G_{\mu }$ which can be
used to define a symplectic structure on $J^{-1}(O)/G$, \cite%
{cushman-sniatycki}.

Consider next the product $P\times O$ with projections $pr_{1}:P\times
O\rightarrow P$ and $pr_{2}:P\times O\rightarrow O$. The form 
\[
\omega _{P\times O}=pr_{1}^{\ast }\omega -pr_{2}^{\ast }\omega _{O}
\]%
on $P\times O$ is symplectic. The Hamiltonian action of $G$ on $P$ with a
momentum map $J:P\rightarrow \mathfrak{g}^{\ast }$ induces a Hamiltonian
action 
\[
G\times (P\times O)\rightarrow P\times O:(g,(p,\nu ))\mapsto (\Phi _{g}(p),{%
Ad}_{g^{-1}}^{\ast }\nu ),
\]%
with a momentum map 
\begin{equation}
J_{P\times O}:P\times O\rightarrow \mathfrak{g}^{\ast }:(p,\nu )\mapsto
J(p)-\nu .  \label{JPO}
\end{equation}%
There is a natural homomorphism between $J^{-1}(O)$ and $J_{P\times
O}^{-1}(0)$ intertwining the actions of $G$ on both spaces. If the action of 
$G$ on $P$ is free and proper, then the orbit spaces $J^{-1}(\mu )/G_{\mu }$
and $J_{P\times O}^{-1}(0)/G$ are symplectomorphic, \cite%
{guillemin-sternberg84}. For a free and proper action, the isomorphic
symplectic manifolds $J^{-1}(\mu )/G_{\mu }$, $J^{-1}(O)/G$ and $J_{P\times
O}^{-1}(0)/G$ are interpreted as reduced phase spaces.

Algebraic reduction at a non-zero orbit $O$ was originally defined by Kimura 
\cite{kimura93} and Wilbour \cite{wilbour93} as follows. Consider the ideal $%
\mathcal{J}^{\mu }$ in $C^{\infty }(P)$ generated by components of $J-\mu
:P\rightarrow \mathfrak{g}^{\ast }.$ In other words,%
\[
\mathcal{J}^{\mu }=\left\{ \sum_{j=1}^{k}f_{j}\langle J-\mu \mid \xi
_{j}\rangle \mid f_{1},...,f_{k}\in C^{\infty }(P)\right\} . 
\]%
Let 
\[
\mathcal{N(J}^{\mu }\mathcal{)}=\{f\in C^{\infty }(P)\mid \{f,h\}\in 
\mathcal{J}^{\mu }\text{ for all }h\in \mathcal{J}^{\mu }\}. 
\]%
The quotient $\mathcal{N(J}^{\mu }\mathcal{)}/(\mathcal{J}^{\mu }\cap 
\mathcal{N(J}^{\mu }\mathcal{))}$ has the structure of a Poisson algebra, It
is interpreted as the Poisson algebra obtained by algebraic reduction of $%
J^{-1}(\mu )$.

Using analogy with regular reduction, Arms defined algebraic reduction at $O$
to be algebraic reduction at $0\in \mathfrak{g}^{\ast }$ of the Hamiltonian
action of $G$ on $P\times O$. If $\mathcal{J}_{P\times O}$ denotes the ideal
generated by components of $J_{P\times O}$, algebraic reduction of $%
J_{P\times O}^{-1}(0)$ gives a Poisson algebra $(C^{\infty }(P\times O)/%
\mathcal{J}_{P\times O})^{G}$. Arms proved that if $G_{\mu }$\ is connected,
then Poisson algebras $\mathcal{N(J}^{\mu }\mathcal{)}/(\mathcal{J}^{\mu
}\cap \mathcal{N(J}^{\mu }\mathcal{))}$ and $(C^{\infty }(P\times O)/%
\mathcal{J}_{P\times O})^{G}$ are isomorphic, \cite{arms}.

In the following we use Arms' definition of algebraic reduction at non-zero
co-adjoint orbits. In other words, for a non-zero co-adjoint orbit $O$, we
interpret $(C^{\infty }(P\times O)/\mathcal{J}_{P\times O})^{G}$ as the
Poisson algebra obtained by algebraic reduction at $O$.

\section{Reduced quantization at $0\in \mathfrak{g}^{\ast }$}

In this section we provide an algebraic analogue of geometric quantization
of reduced Poisson algebras described in Section 2.4. The basic idea is to
replace the prequantization line bundle $\pi :L\rightarrow P$ be the space $%
S^{\infty }(L)$ of sections of $\pi $ and use techniques of algebraic
geometry to construct the reduced space and the reduced quantization map.

\subsection{Prequantization}

We consider here a prequantization structure on $(P,\omega )$ given by a
complex line bundle $\pi :L\rightarrow P$, with connection $\nabla $ such
that the curvature form of $\nabla $ is $-h^{-1}\omega $, where $h$ is
Planck's constant, and a connection invariant Hermitian form $\langle \cdot
|\cdot \rangle $. The corresponding \emph{prequantization} map is given by 
\begin{equation}
\mathcal{P}:C^{\infty }(P)\times S^{\infty }(L)\rightarrow S^{\infty
}(L):(f,\sigma )\mapsto \mathcal{P}_{f}\sigma =(-i\hbar \nabla
_{X_{f}}+f)\sigma ,  \label{10}
\end{equation}%
where $\hbar =h/2\pi $. For each $f\in C^{\infty }(P)$, such that the
Hamiltonian vector field $X_{f}$ is complete, the operator $\mathcal{P}_{f}$
is skew adjoint on the space of sections of $\lambda $ that are square
integrable with respect to the scalar product 
\begin{equation}
\left( \sigma _{1}|\sigma _{2}\right) =\int_{P}\langle \sigma _{1}|\sigma
_{2}\rangle \omega ^{n},  \label{ScalarProduct}
\end{equation}%
where $n=\frac{1}{2}\dim P$.

Restricting the prequantization map $\mathcal{P}$ to the Poisson algebra
spanned by $J_{\xi }$, for $\xi \in \mathfrak{g}$, we get a representation $%
\xi \mapsto (-i\hbar )^{-1}\mathcal{P}_{J_{\xi }}$ of $\mathfrak{g}$ on $%
S^{\infty }(L).$ This representation integrates to a representation%
\[
\mathcal{T}:G\times S^{\infty }(L)\rightarrow S^{\infty }(L):(g,\sigma
)\mapsto \mathcal{T}_{g}\sigma 
\]%
of $G$ on $S^{\infty }(L),$ called the \emph{prequantization representation}%
. The prequantization representation $\mathcal{T}$ is unitary on the space
of sections of $\pi :L\rightarrow P$ which are square integrable with
respect to the scalar product (\ref{ScalarProduct}).

The space $S^{\infty }(L)$ of sections of the line bundle is a module over
the ring $C^{\infty }(P)$ of smooth functions on $P$. We want to construct
from $S^{\infty }(L)$ a module over the quotient ring $C^{\infty }(P)/%
\mathcal{J}$. This construction is an algebraic analogue of the restriction
of the line bundle to $J^{-1}(0)$. Let 
\[
\mathcal{J}S^{\infty }(L)=\mathrm{span}\{f\sigma \mid f\in \mathcal{J}\text{%
, }\sigma \in S^{\infty }(L)\}. 
\]%
For every $[\sigma ]\in S^{\infty }(L)/\mathcal{J}S^{\infty }(L)$ and $%
[f]\in C^{\infty }(P)/\mathcal{J}$, the class $[f\sigma ]$ in $S^{\infty
}(L)/\mathcal{J}S^{\infty }(L)$ is independent of representatives $f$ of $%
[f] $ and $\sigma $ of $[\sigma ]$. Hence, $S^{\infty }(L)/\mathcal{J}%
S^{\infty }(L)$ is a module over $C^{\infty }(P)/\mathcal{J}$ with $%
[f][\sigma ]=[f\sigma ].$

\begin{proposition}
\label{proposition1}The prequantization representation $\mathcal{T}$ of $G$
on $S^{\infty }(L)$ induces a representation of $G$ on $S^{\infty }(L)/%
\mathcal{J}S^{\infty }(L)$ given by 
\begin{equation}
G\times S^{\infty }(L)/\mathcal{J}S^{\infty }(L)\rightarrow S^{\infty }(L)/%
\mathcal{J}S^{\infty }(L):(g,[\sigma ])\mapsto \lbrack \mathcal{T}_{g}\sigma
].  \label{12}
\end{equation}%
A class $[\sigma ]\in S^{\infty }(L)/\mathcal{J}S^{\infty }(L)$ is $G$%
-invariant if 
\[
\mathcal{P}_{J_{\xi }}\sigma \in \mathcal{J}S^{\infty }(L)\text{ for all }%
\xi \in \mathfrak{g}\text{.}
\]
\end{proposition}

\begin{proof}
In order to show that the map in equation (\ref{12}) is well defined we need
to show that, for every $g\in G$, the map $\mathcal{T}_{g}:S^{\infty
}(L)\rightarrow S^{\infty }(L)$ preserves $\mathcal{J}S^{\infty }(L)$. Since 
$G$ is connected, it suffices to show that $\mathcal{P}_{J_{\xi }}\sigma \in 
\mathcal{J}S^{\infty }(L)$ for all $\sigma \in \mathcal{J}S^{\infty }(L)$
and $\xi \in \mathfrak{g}$. Consider $J_{\zeta }\sigma \in \mathcal{J}%
S^{\infty }(L)$. Using equation (\ref{10}) we get, for each $\xi \in 
\mathfrak{g}$, 
\begin{eqnarray*}
\mathcal{P}_{J_{\xi }}(J_{\zeta }\sigma ) &=&(-i\hbar \nabla _{X_{J_{\xi
}}}+J_{\xi })(J_{\zeta }\sigma )=J_{\zeta }(-i\hbar \nabla _{X_{J_{\xi
}}}+J_{\xi })\sigma -i\hbar (X_{J_{\xi }}J_{\zeta })\sigma \\
&=&J_{\zeta }(-i\hbar \nabla _{X_{J_{\xi }}}+J_{\xi })\sigma -i\hbar J_{[\xi
,\zeta ]}\sigma \in \mathcal{J}S^{\infty }(L).
\end{eqnarray*}%
Hence, $\mathcal{T}_{g}$ maps $\mathcal{J}S^{\infty }(L)$ to itself for all $%
g\in G$. Thus, the map $(g,[\sigma ])\mapsto \lbrack \mathcal{T}_{g}\sigma ]$
is well defined for all $g\in G$. It is a representation of $G$ on $%
S^{\infty }(L)/\mathcal{J}S^{\infty }(L)$ because $\mathcal{T}$ is a
representation.

Since $G$ is connected, $G$-invariance of $[\sigma ]\in S^{\infty }(L)/%
\mathcal{J}S^{\infty }(L)$ is a consequence of invariance of $[\sigma ]$
under one parameter sugroups of $G$. Hence, \thinspace $\lbrack \sigma ]$ is 
$G$-invariant if $[\mathcal{P}_{J_{\xi }}\sigma ]=0$ for all $\xi \in 
\mathfrak{g}$. This is equivalent to $\mathcal{P}_{J_{\xi }}\sigma \in 
\mathcal{J}S^{\infty }(L)$ for all $\xi \in \mathfrak{g}$.
\end{proof}

In the case of prequantization of regular reduction at $0\in \mathfrak{g}%
^{\ast }$, the reduced space of wave functions is the space of $G$-invariant
sections The next step of prequantization corresponds to pushing forward $%
L_{\mid J^{-1}(0)}$ to $J^{-1}(0)/G$. It is usually done in terms of $G$%
-invariant sections of $L_{\mid J^{-1}(0)}$. An algebraic reduction analogue
of the space of $G$-invariant sections of $L_{\mid J^{-1}(0)}$ is the space%
\begin{equation}
(S^{\infty }(L)/\mathcal{J}S^{\infty }(L))^{G}=\{[\sigma ]\in S^{\infty }(L)/%
\mathcal{J}S^{\infty }(L)\mid \mathcal{P}_{J_{\xi }}\sigma \in JS^{\infty
}(L)\}  \label{13}
\end{equation}%
of $G$-invariant elements of $S^{\infty }(L)/\mathcal{J}S^{\infty }(L).$ It
plays the role of the space of wave functions in quantization of algebraic
reduction at $0\in \mathfrak{g}^{\ast }$.

\begin{theorem}
\label{theorem1}$(S^{\infty }(L)/\mathcal{J}S^{\infty }(L))^{G}$ is a module
over the ring $(C^{\infty }(P)/\mathcal{J})^{G}$. Moreover, the \emph{%
algebraically reduced prequantization map} 
\[
(C^{\infty }(P)/\mathcal{J)}^{G}\mathcal{\times }{\left( S^{\infty }(L)/%
\mathcal{J}S^{\infty }(L)\right) }^{G}\rightarrow {\left( S^{\infty }(L)/%
\mathcal{J}S^{\infty }(L)\right) }^{G}:([f],[\sigma ])\mapsto \mathcal{P}%
_{[f]}[\sigma ] 
\]%
given by 
\begin{equation}
\mathcal{P}_{[f]}[\sigma ]=[\mathcal{P}_{f}\sigma ],  \label{13.1}
\end{equation}%
where $f$ is any representative of $[f]$ and $\sigma $ is any representative
of $[\sigma ]$, is well defined.
\end{theorem}

\begin{proof}
We know that $S^{\infty }(L)/\mathcal{J}S^{\infty }(L)$ is a module over $%
C^{\infty }(P)/\mathcal{J}.$ If $[f]\in (C^{\infty }(P)/\mathcal{J})^{G}$
and $[\sigma ]\in (S^{\infty }(L)/\mathcal{J}S^{\infty }(L))^{G}$ then $%
X_{J_{\xi }}f\in \mathcal{J}$ and $\mathcal{P}_{J_{\xi }}\sigma \in \mathcal{%
J}S^{\infty }(L)$ for all $\xi \in \mathfrak{g}$. Hence, $\mathcal{P}%
_{J_{\xi }}(f\sigma )=-i\hbar (X_{J_{\xi }}f)\sigma +f\mathcal{P}_{J_{\xi
}}\sigma \in \mathcal{J}\Gamma ^{\infty }(L)$ for all $\xi \in \mathfrak{g}$%
. This implies that $[f][\sigma ]=[f\sigma ]$ is $G$-invariant. Hence, $%
(S^{\infty }(L)/\mathcal{J}S^{\infty }(L))^{G}$ is a module over $(C^{\infty
}(P)/\mathcal{J})^{G}$.

To show that the reduced prequantization map is well defined it suffices to
show that, for every $[\sigma ]\in \left( S^{\infty }(L)/\mathcal{J}%
S^{\infty }(L)\right) ^{G}$ and $[f]\in \left( C^{\infty }(P)/\mathcal{J}%
\right) ^{G}$, the equivalence class $[\mathcal{P}_{f}\sigma ]$ is in $%
\left( S^{\infty }(L)/\mathcal{J}S^{\infty }(L)\right) ^{G}$ and it does not
depend on the choice of representatives $f$ of $[f]$ and $\sigma $ of $%
[\sigma ].$ Since $[\sigma ]\in \left( S^{\infty }(L)/\mathcal{J}S^{\infty
}(L)\right) ^{G}$ it follows that, for every representative $\sigma $ of $%
[\sigma ]$, $\mathcal{P}_{J_{\xi }}\sigma \in \mathcal{J}S^{\infty }(L)$ for
all $\xi \in \mathfrak{g}$. Similarly, $[f]\in \left( C^{\infty }(P)/%
\mathcal{J}\right) ^{G}$ implies that $X_{J_{\xi }}f\in \mathcal{J}$ for
every representative $f$ of $[f]$. For $[f]\in \left( C^{\infty }(P)/%
\mathcal{J}\right) ^{G}$ and $J_{\xi }\sigma \in \mathcal{J}S^{\infty }(L),$ 
\begin{eqnarray*}
\mathcal{P}_{f}(J_{\xi }\sigma ) &=&(-i\hbar \nabla _{X_{f}}+f)(J_{\xi
}\sigma )=J_{\xi }(-i\hbar \nabla _{X_{f}}+f)\sigma -i\hbar (X_{f}J_{\xi
})\sigma  \\
&=&J_{\xi }(-i\hbar \nabla _{X_{f}}+f)\sigma +i\hbar (X_{J_{\xi }}f)\sigma
\in \mathcal{J}\Gamma ^{\infty }(L).
\end{eqnarray*}%
This implies that, for $f\in \lbrack f]$, $\mathcal{P}_{f}$ maps $\mathcal{J}%
S^{\infty }(L)$ to itself. Hence, $[\mathcal{P}_{f}\sigma ]$ does not depend
the representative $\sigma $ of $[\sigma ]$. For $J_{\xi }f\in \mathcal{J}$,
and $[\sigma ]\in \left( S^{\infty }(L)/\mathcal{J}S^{\infty }(L)\right) ^{G}
$, equation (\ref{10}) yields 
\begin{eqnarray*}
\mathcal{P}_{J_{\xi }f}\sigma  &=&(-i\hbar \nabla _{X_{J_{\xi }f}}+J_{\xi
}f)\sigma  \\
&=&-i\hbar J_{\xi }\nabla _{X_{f}}\sigma +i\hbar f\mathcal{P}_{J_{\xi
}}\sigma \in \mathcal{J}S^{\infty }(L).
\end{eqnarray*}%
Hence, $[\mathcal{P}_{f}\sigma ]$ does not depend on the representative $f$
of $[f]$. Combining these results, we conclude that the equivalence class $[%
\mathcal{P}_{f}\sigma ]$ depends only on $[\sigma ]\in \left( S^{\infty }(L)/%
\mathcal{J}S^{\infty }(L)\right) ^{G}$ and $[f]\in \left( C^{\infty }(P)/%
\mathcal{J}\right) ^{G}$.

It remains to show that $[\mathcal{P}_{f}\sigma ]\in (S^{\infty }(L)/%
\mathcal{J}S^{\infty }(L))^{G}$. For, $\xi \in \mathfrak{g,}$ $[f]\in \left(
C^{\infty }(P)/\mathcal{J}\right) ^{G}$ and $[\sigma ]\in \left( S^{\infty
}(L)/\mathcal{J}S^{\infty }(L)\right) ^{G},$%
\[
\mathcal{P}_{J_{\xi }}\mathcal{P}_{f}\sigma =\mathcal{P}_{f}\mathcal{P}%
_{J_{\xi }}\sigma +[\mathcal{P}_{J_{\xi }},\mathcal{P}_{f}]\sigma .
\]%
But $[\sigma ]\in \left( {S}^{\infty }(L)/\mathcal{J}{S}^{\infty }(L)\right)
^{G}$ implies $\mathcal{P}_{J_{\xi }}\sigma \in \mathcal{J}{S}^{\infty }(L)$
so that $\mathcal{P}_{f}\mathcal{P}_{J_{\xi }}\sigma \in \mathcal{J}{S}%
^{\infty }(L)$ by the first part of the proof. On the other hand,%
\[
\lbrack \mathcal{P}_{J_{\xi }},\mathcal{P}_{f}]\sigma =\mathcal{-}i\hbar 
\mathcal{P}_{\{J_{\xi },f\}}\sigma =\mathcal{-}i\hbar (-i\hbar \nabla
_{X_{\{J_{\xi },f\}}}+\{J_{\xi },f\})\sigma .
\]%
Since $[f]\in \left( C^{\infty }(P)/\mathcal{J}\right) ^{G}$ it follows that 
$\{J_{\xi },f\}=-X_{J_{\xi }}f\in \mathcal{J}$. Hence, $\{J_{\xi
},f\}=\sum_{j}f_{j}J_{\zeta _{j}}$ for some $f_{j}\in C^{\infty }(P)$ and $%
\zeta _{j}\in \mathfrak{g}$. Moreover, $%
X_{f_{1}f_{2}}=f_{1}X_{f_{2}}+f_{2}X_{f_{1}}$ implies that 
\[
\sum_{j}\nabla _{X_{f_{j}J_{\zeta _{j}}}}=\sum_{j}(f_{j}\nabla _{X_{J_{\zeta
_{j}}}}+J_{\zeta _{j}}\nabla _{X_{f_{j}}}).
\]%
Therefore,%
\begin{eqnarray*}
(-i\hbar \nabla _{X_{\{J_{\xi },f\}}}+\{J_{\xi },f\})\sigma 
&=&\sum_{j}\left( -i\hbar (f_{j}\nabla _{X_{J_{\zeta _{j}}}}+J_{\zeta
_{j}}\nabla _{X_{f_{j}}})+f_{j}J_{\zeta _{j}}\right) \sigma  \\
&=&\sum_{j}\left( (f_{j}\mathcal{P}_{J_{\zeta _{j}}}-i\hbar J_{\zeta
_{j}}\nabla _{X_{f_{j}}}\right) \sigma \in \mathcal{J}S^{\infty }(L),
\end{eqnarray*}%
because $[\sigma ]\in \left( S(L)/\mathcal{J}S^{\infty }(L)\right) ^{G}$.
Hence, $[\mathcal{P}_{f}\sigma ]\in \left( S^{\infty }(L)/\mathcal{J}%
S^{\infty }(L)\right) ^{G}$. This completes the proof of Theorem \ref%
{theorem1}.
\end{proof}

\subsection{Polarization}

A \ \emph{polarization} of a symplectic manifold $(P,\omega )$ is an
involutive complex Lagrangian distribution $F$ on $P$. Let $\bar{F}$ denote
the complex conjugate of $F$. A polarization $F$ is \emph{strongly admissible%
} if $D=F\cap \bar{F}\cap TP$ and $E=(F+\bar{F})\cap TP$ are involutive
distributions, spaces $P/D$ and $P/E$ of leaves of $D$ and $E$,
respectively, are quotient manifolds of $P$ and the canonical projection $%
X/D\rightarrow X/E$ is a submersion. For a strongly admissible polarization $%
F$, leaves of $D$ have structure of affine manifolds. If leaves of $D$ are
complete in their affine structure, we say that $F$ is \emph{complete}.

We assume here that $(P,\omega )$ is endowed with a strongly admissible
complete polarization $F$. The space of \emph{polarized sections} of $L$ is 
\[
S_{F}^{\infty }(L)=\{\sigma \in S^{\infty }(L)\mid \nabla _{u}\sigma =0\text{
for all }u\in F\}.
\]%
It is a module over the ring $C_{F}^{\infty }(P)^{0}$ consisting of smooth
functions on $P$ which are annihilated by $F$. In other words, 
\[
C_{F}^{\infty }(P)^{0}=\{f\in C^{\infty }(P)\mid uf=0\text{ for all }u\in
F\}.
\]

Functions in $C^{\infty }(P)$, such that their Hamiltonian vector fields
preserve $F$, form a Poisson subalgebra $C_{F}^{\infty }(P)$ of $C^{\infty
}(P).$ In other words, a function $f$ is in $C_{F}^{\infty }(P)$ if, for
every vector field $X$ on $P$ with values in $F$, the bracket $[X_{f},X]$
has values in $F$. Clearly, $C_{F}^{\infty }(P)^{0}$ is contained in $%
C_{F}^{\infty }(P).$ A function $f\in C^{\infty }(P)$ is \emph{quantizabl}e
in terms of the polarization $F$ if it belongs to $C_{F}^{\infty }(P)$. The 
\emph{quantization }map $\mathcal{Q}$ is obtained from the prequantization
map $\mathcal{P}$ by restricting its domain to $C_{F}^{\infty }(P)\times
S_{F}^{\infty }(L)$ and co-domain to $S_{F}^{\infty }(L)$. In other words, 
\begin{equation}
\mathcal{Q}:C_{F}^{\infty }(P)\times S_{F}^{\infty }(L)\rightarrow
S_{F}^{\infty }(L):(f,\sigma )\mapsto \mathcal{Q}_{f}\sigma =\mathcal{P}%
_{f}\sigma .  \label{14}
\end{equation}

We assume that the action $\Phi $ of $G$ on $P,$ given by (\ref{4}),
preserves $F$. Then, for every $\xi \in \mathfrak{g}$, the momentum $J_{\xi
} $ is in $C_{F}^{\infty }(P)$ and the map $\xi \mapsto (-i\hbar )^{-1}%
\mathcal{Q}_{J_{\xi }}$ is a representation of $\mathfrak{g}$ on $%
S_{F}^{\infty }(L)$. This representation integrates to a linear
representation of $G$ on $S_{F}^{\infty }(L)$ denoted by 
\begin{equation}
\mathcal{R}:G\times S_{F}^{\infty }(L)\rightarrow S_{F}^{\infty
}(L):(g,\sigma )\mapsto \mathcal{R}_{g}\sigma ,  \label{R}
\end{equation}%
and called a \emph{quantization} representation of $G$. Since the group $G$
is connected, a section $\sigma \in S_{F}^{\infty }(L)$ is $G$-invariant if $%
\mathcal{Q}_{J_{\xi }}\sigma =0$ for all $\xi \in \mathfrak{g}$.

Let $\Sigma :S^{\infty }(L)\rightarrow S^{\infty }(L)/\mathcal{J}S^{\infty
}(L):\sigma \mapsto \Sigma (\sigma )=[\sigma ]$ be the natural projection.
We denote by $S_{F}^{\infty }(L)/\mathcal{J}S^{\infty }(L)$ the image of $%
S_{F}^{\infty }(L)$ under $\Sigma $. It can be identified with $%
S_{F}^{\infty }(L)/(S_{F}^{\infty }(L)\cap \mathcal{J}S^{\infty }(L))$. In
other words, 
\[
S_{F}^{\infty }(L)/\mathcal{J}S^{\infty }(L)=\Sigma (S_{F}^{\infty
}(L))=S_{F}^{\infty }(L)/(S_{F}^{\infty }(L)\cap \mathcal{J}S_{F}^{\infty
}(L)). 
\]%
The representation space of quantization of algebraic reduction at $0\in 
\mathfrak{g}^{\ast }$ is the space 
\begin{equation}
(S_{F}^{\infty }(L)/\mathcal{J}S^{\infty }(L))^{G}=\{[\sigma ]\in
S_{F}^{\infty }(L)/\mathcal{J}S^{\infty }(L)\mid \mathcal{Q}_{J_{\xi
}}\sigma \in \mathcal{J}S^{\infty }(L)\text{ for all }\xi \in \mathfrak{g\}}
\label{SFLG}
\end{equation}%
of $G$-invariant elements in $S_{F}^{\infty }(L)/\mathcal{J}S^{\infty }(L)$.

Let $C_{F}^{\infty }(P)^{0}/\mathcal{J}$ be the image of $C_{F}^{\infty
}(P)^{0}$ under the natural projection $C^{\infty }(P)\rightarrow C^{\infty
}(P)/\mathcal{J}$. We denote by $(C_{F}^{\infty }(P)^{0}/\mathcal{J})^{G}$
the space of invariant elements in $C_{F}^{\infty }(P)^{0}/\mathcal{J}$.

\begin{proposition}
\label{proposition2}$(C_{F}^{\infty }(P)^{0}/\mathcal{J})^{G}$ is a ring and 
$(S_{F}^{\infty }(L)/\mathcal{J}S^{\infty }(L))^{G}$ is a module over $%
(C_{F}^{\infty }(P)^{0}/\mathcal{J})^{G}$.
\end{proposition}

\begin{proof}
Ring structure of $(C_{F}^{\infty }(P)^{0}/\mathcal{J})^{G}$ follows from
the fact that $C_{F}^{\infty }(P)^{0}$ is a ring. By definition, $[f]\in
(C_{F}^{\infty }(P)^{0}/\mathcal{J})^{G\text{ }}$ has a representative $f\in
C_{F}^{\infty }(P)$ such that $X_{J_{\xi }}f\in \mathcal{J}$ for all $\xi
\in \mathfrak{g}$. Similarly, $[\sigma ]\in {\left( S_{F}^{\infty }(L)/%
\mathcal{J}S^{\infty }(L)\right) }^{G}$ has a representative in $%
S_{F}^{\infty }(L)$ such that $\mathcal{P}_{J_{\xi }}\sigma \in \mathcal{J}%
S^{\infty }(L)$ for all $\xi \in \mathfrak{g}.$ Since $[f][\sigma ]=[f\sigma
],$ where $f\sigma \in S_{F}^{\infty }(L)$ and $\mathcal{P}_{J_{\xi
}}(f\sigma )=f\mathcal{P}_{J_{\xi }}(\sigma )-i\hbar (X_{\xi }f)\sigma \in 
\mathcal{J}S^{\infty }(L)$ for all $\xi \in \mathfrak{g}$, it follows that $%
[f][\sigma ]\in {\left( S_{F}^{\infty }(L)/\mathcal{J}S^{\infty }(L)\right) }%
^{G}$. Therefore, ${\left( S_{F}^{\infty }(L)/\mathcal{J}S^{\infty
}(L)\right) }^{G}$ is a module over $(C_{F}^{\infty }(P)^{0}/\mathcal{J})^{G%
\text{ }}$.
\end{proof}

An element $[f]$ of $(C^{\infty }(P)/J)^{G}$ is \emph{quantizable} in terms
of the polarization $F$ if it has a representative $f\in C_{F}^{\infty }(P)$%
. Quantizable\emph{\ }elements of $(C^{\infty }(P)/\mathcal{J})^{G}$ form a
Poisson subalgebra 
\[
(C_{F}^{\infty }(P)/\mathcal{J})^{G}=\{[f]\in (C^{\infty }(P)/\mathcal{J}%
)^{G}\mid f\in C_{F}^{\infty }(P)\}.
\]

\begin{theorem}
\label{theorem2}The \emph{algebraically reduced quantization} map 
\[
(C_{F}^{\infty }(P)/\mathcal{J)}^{G}\mathcal{\times }{\left( S_{F}^{\infty
}(L)/\mathcal{J}S^{\infty }(L)\right) }^{G}\rightarrow {\left( S_{F}^{\infty
}(L)/\mathcal{J}S^{\infty }(L)\right) }^{G}:([f],[\sigma ])\mapsto \mathcal{Q%
}_{[f]}[\sigma ] 
\]%
given by 
\[
\mathcal{Q}_{[f]}[\sigma ]=[\mathcal{Q}_{f}\sigma ] 
\]%
is well defined. Moreover, 
\[
\mathcal{Q}_{[f]}[\sigma ]=\mathcal{P}_{[f]}[\sigma ]=[\mathcal{P}_{f}\sigma
], 
\]%
for all $[f]\in (C_{F}^{\infty }(P)/\mathcal{J)}^{G}$ and $[\sigma ]\in {%
\left( S_{F}^{\infty }(L)/\mathcal{J}S^{\infty }(L)\right) }^{G}$.
\end{theorem}

\begin{proof}
By definition, $[f]\in (C_{F}^{\infty }(P)/\mathcal{J})^{G\text{ }}$ has a
representative $f\in C_{F}^{\infty }(P)$ and the Hamiltonian vector field $%
X_{f}$ of $f$ preserves $F$. Similarly, $[\sigma ]\in \left( S_{F}^{\infty
}(L)/\mathcal{J}S^{\infty }(L)\right) ^{G}$ has a representative $\sigma \in 
$ $S_{F}^{\infty }(L)$. It follows that $\mathcal{Q}_{f}\sigma $ is in $%
S_{F}^{\infty }(L)$, so that $[\mathcal{Q}_{f}\sigma ]\in S_{F}^{\infty }(L)/%
\mathcal{J}S^{\infty }(L)$. But, $\mathcal{Q}_{f}\sigma =\mathcal{P}%
_{f}\sigma $, and we have shown in Theorem \ref{theorem1} that $[\mathcal{P}%
_{f}\sigma ]\in \left( S^{\infty }(L)/\mathcal{J}S^{\infty }(L)\right) ^{G}$
and that $[\mathcal{P}_{f}\sigma ]$ is independent of the choice of
representatives of $[f]$ and of $[\sigma ]$. Hence, $[\mathcal{Q}_{f}\sigma ]
$ lies in 
\[
\left( S^{\infty }(L)/\mathcal{J}S^{\infty }(L)\right) ^{G}\cap
(S_{F}^{\infty }(L)/\mathcal{J}S^{\infty }(L))=({S}_{F}^{\infty }(L)/%
\mathcal{J}S^{\infty }(L))^{G}
\]%
and it is independent of the choice of representative of $[f]$ and of $%
[\sigma ]$. Thus, the reduced quantization map is well defined. On the other
hand, $\mathcal{Q}_{f}\sigma =\mathcal{P}_{f}\sigma $ implies that%
\[
\mathcal{Q}_{[f]}[\sigma ]=[\mathcal{Q}_{f}\sigma ]=[\mathcal{P}_{f}\sigma ]=%
\mathcal{P}_{[f]}[\sigma ].
\]%
This completes proof of Theorem \ref{theorem2}.
\end{proof}

\subsection{K\"{a}hler polarizations}

A polarization $F$ is called \emph{K\"{a}hler} if it is purely imaginary,
i.e. $\bar{F}\cap F=0$, and $i\omega (w,\bar{w})>0$ for all non-zero $w\in F$%
.\ In this case, $(P,\omega )$ is a K\"{a}hler manifold, $F$ is the
distribution of antiholomorphic directions, $L$ is a holomorphic line bundle
over $P$ and $S_{F}^{\infty }(L)$ consists of holomorphic sections.
Holomorphic sections of $L$ which are square integrable with respect to the
volume form $\omega ^{n}$, where $2n=\dim P$, form a Hilbert space $\mathcal{%
H}$. In this case, the linear representation $\mathcal{R}$ of $G$ on $%
S_{F}^{\infty }(L)$ restricts to a unitary representation $\mathcal{U}$ of $G
$ on $\mathcal{H}$.

\begin{lemma}
\label{lemma1}Let $F$ be a K\"{a}hler polarization of $(P,\omega )$ and let $%
J:P\rightarrow \mathfrak{g}^{\ast }$ be a momentum map for a Hamiltonian
action of $G$ on $P$. If $J^{-1}(0)$ contains a Lagrangian submanifold of $%
(P,\omega )$ then $S_{F}^{\infty }(L)\cap \mathcal{J}S^{\infty }(L)=0.$
\end{lemma}

\begin{proof}
Each $\sigma \in S_{F}^{\infty }(L)\cap \mathcal{J}S^{\infty }(L)$ is a
holomorphic section of $L$ which vanishes on $J^{-1}(0)$. Since $J^{-1}(0)$
contains a Lagrangian submanifold it follows that all derivatives of $\sigma 
$ vanish on $J^{-1}(0)$. Hence, $\sigma =0$. 
\hfill%
\end{proof}

Let $S_{F}^{\infty }(L)^{G}$ denote the space of $G$-invariant sections of $%
S_{F}^{\infty }(L)$. Since $G$ is connected, it follows that 
\[
S_{F}^{\infty }(L)^{G}=\{\sigma \in S_{F}^{\infty }(L)\mid \mathcal{Q}%
_{J_{\xi }}\sigma =0\text{ for all }\xi \in \mathfrak{g}\}. 
\]

\begin{theorem}
\label{theorem3}Let $F$ be a $G$-invariant K\"{a}hler polarization of $%
(P,\omega )$ and let $J:P\rightarrow \mathfrak{g}^{\ast }$ be a momentum map
for the action of $G$ on $(P,\omega )$ such that $J^{-1}(0)$ contains a
Lagrangian submanifold of $(P,\omega ).$ Then, the map%
\[
\Xi :S_{F}^{\infty }(L)^{G}\rightarrow \left( S_{F}^{\infty }(L)/\mathcal{J}%
S^{\infty }(L)\right) ^{G}:\sigma \mapsto \lbrack \sigma ] 
\]%
is an isomorphism of vector spaces such that, for each $G$-invariant
function $f\in C_{F}^{\infty }(P),$ 
\[
\Xi \circ \mathcal{Q}_{f}=\mathcal{Q}_{[f]}\circ \Xi . 
\]
\end{theorem}

\begin{proof}
By definition, $S_{F}^{\infty }(L)/\mathcal{J}S^{\infty }(L)=\Sigma
(S_{F}^{\infty }(L))$, where $\Sigma $ is the natural projection from $%
S^{\infty }(L)$ to $S^{\infty }(L)/\mathcal{J}S^{\infty }(L)$. If $\sigma
\in S_{F}^{\infty }(L)$ is $G$-invariant, then $\mathcal{Q}_{J_{\xi }}\sigma
=0\in \mathcal{J}S^{\infty }(L)$ for all $\xi \in \mathfrak{g}$ and equation
(\ref{SFLG}) implies that $[\sigma ]\in \left( S_{F}^{\infty }(L)/\mathcal{J}%
S^{\infty }(L)\right) ^{G}$. Hence, $\Xi $ is the restriction of $\Sigma $
to domain $S_{F}^{\infty }(L)^{G}$ and codomain $\left( S_{F}^{\infty }(L)/%
\mathcal{J}S^{\infty }(L)\right) ^{G}.$

We first show that $\Xi $ is one-to one. Suppose that $[\sigma ]\in
S_{F}^{\infty }(L)/\mathcal{J}S_{F}^{\infty }(L)$ has two representatives $%
\sigma $ and $\sigma ^{\prime }$ in $S_{F}^{\infty }(L)$. Since $%
S_{F}^{\infty }(L)$ is a vector space, it follows that $\sigma -\sigma
^{\prime }$ is in $S_{F}^{\infty }(L)$. On the other hand, $[\sigma
]=[\sigma ^{\prime }]$ implies that $\sigma -\sigma ^{\prime }\in \mathcal{J}%
S^{\infty }(L)$. Hence, $\sigma -\sigma ^{\prime }\in S_{F}^{\infty }(L)\cap 
\mathcal{J}S^{\infty }(L)$. Lemma \ref{lemma1} implies that $\sigma -\sigma
^{\prime }=0$. Thus each element $[\sigma ]\in \left( S_{F}^{\infty }(L)/%
\mathcal{J}S^{\infty }(L)\right) $ has a unique representative $\sigma \in
S_{F}^{\infty }(L)$.

We need to show that $\Xi $ is onto $\left( S_{F}^{\infty }(L)/\mathcal{J}%
S^{\infty }(L)\right) ^{G}$. Each $[\sigma ]\in \left( S_{F}^{\infty }(L)/%
\mathcal{J}S^{\infty }(L)\right) ^{G}$ has a representative $\sigma \in
S_{F}^{\infty }(L)$ such that $\mathcal{Q}_{J_{\xi }}\sigma \in \mathcal{J}%
S^{\infty }(L)$ for all $\xi \in \mathfrak{g}$. On the other hand, $\mathcal{%
Q}_{J_{\xi }}\sigma \in S_{F}^{\infty }(L)$. Hence, $\mathcal{Q}_{J_{\xi
}}\sigma \in S_{F}^{\infty }(L)\cap \mathcal{J}S^{\infty }(L)=0.$ Thus, $%
\sigma $ is $G$-invariant which implies that $\Xi $ maps $S_{F}^{\infty
}(L)^{G}$ onto $\left( S_{F}^{\infty }(L)/\mathcal{J}S^{\infty }(L)\right)
^{G}$.

For every $G$-invariant function $f\in C_{F}^{\infty }(P)$, the class $[f]$
is in $(C_{F}^{\infty }(P)/\mathcal{J})^{G}$. Moreover, for each $\sigma \in
S_{F}^{\infty }(L)^{G}$, $[\sigma ]=\Xi (\sigma )\in \left( S_{F}^{\infty
}(L)/\mathcal{J}S^{\infty }(L)\right) ^{G}$ and $\mathcal{Q}_{f}\sigma \in
S_{F}^{\infty }(L)^{G}.$ Hence, 
\[
\mathcal{Q}_{[f]}\circ \Xi (\sigma )=\mathcal{Q}_{[f]}[\sigma ]=[\mathcal{Q}%
_{f}\sigma ]=\Xi \circ \mathcal{Q}_{f}\sigma ,
\]%
which completes the proof. 
\hfill%
\end{proof}

\subsection{Generalized invariant sections}

In this section we discuss another special case when the polarization $F$ is
real, that is $F=D\otimes \mathbb{C}$, where $D$ is an involutive Lagrangian
distribution on $(P,\omega ).$ We assume, that the space $Q$ of integral
manifolds of $D$ is a manifold of $P$ and the natural projection map $%
\vartheta :P\rightarrow Q$, associating to each $p\in P$ the integral
manifold $\Lambda $ of $D$ through $p$, is a submersion. This implies that,
for each integral manifold $\Lambda $ of $D$, the connection $\nabla $ on $L$
induces a flat connection on the restriction $L_{\mid \Lambda }$ of $L$ to $%
\Lambda $. Further, we assume vanishing of the holonomy of the flat
connection on $L_{\mid \Lambda }$ for every integral manifold $\Lambda $ of $%
D$. Under these assumptions, the parallel transport in $L$ along integral
manifolds of $D$ defines an equivalence relation $\symbol{126}$ on $L$. The
quotient $\tilde{L}=L/\symbol{126}$ has the structure of a complex line
bundle over $Q$. Let $\tilde{\pi}:\tilde{L}\rightarrow Q$ denote the
projection map associating to each equivalence class in $\tilde{L}$ the
corresponding integral manifold $\Lambda \in Q$. Further, let $\lambda
:L\rightarrow \tilde{L}$ be the canonical projection associating to each $%
z\in L$ its equivalence class $\lambda (z)\in \tilde{L}$. Then, $\tilde{\pi}%
\circ \lambda =\vartheta \circ \pi $ and, for each $p\in P$, the map $%
\lambda $ restricts to an isomorphism of the fibre $L_{p}$ onto the fibre $%
\tilde{L}_{\vartheta (p)}.$

Concider a section $\sigma \in S_{F}^{\infty }(L)$. Since $F=D\otimes 
\mathbb{C}$, it follows that $\sigma $ is covariantly constant along $D$,
and $\lambda (\sigma (p_{1}))=\lambda (\sigma (p_{2}))$ if $p_{1}$ and $p_{2}
$ are in the same integral manifold $\Lambda $ of $D$. Hence,  $\sigma \in
S_{F}^{\infty }(L)$ induces a smooth section $\tilde{\sigma}$ of $\tilde{\pi}%
:\tilde{L}\rightarrow Q$ such that $\lambda \circ \sigma =\tilde{\sigma}%
\circ \vartheta $. Conversely, for a section $\tilde{\sigma}$ of $\tilde{\pi}%
:\tilde{L}\rightarrow Q$, there exists a unique section $\sigma $ of $\pi
:L\rightarrow P$ covariantly constant along $D$ and such that $\lambda \circ
\sigma =\tilde{\sigma}\circ \vartheta $. Since $\sigma $ is covariantly
constant along $D$ it is covariantly constant along $F=D\otimes \mathbb{C}$.
We refer to $\sigma $ as the horizontal lift of $\tilde{\sigma}$ and write $%
\sigma =\mathrm{lift}$~$\tilde{\sigma}.$ The map $\tilde{\sigma}\mapsto 
\mathrm{lift~}\tilde{\sigma}$ is a vector space isomorphism of the space $%
S^{\infty }(\tilde{L})$ of smooth sections of $\tilde{L}$ onto $%
S_{F}^{\infty }(L)$. Moreover, $\mathrm{lift~}\tilde{f}\tilde{\sigma}%
=(\vartheta ^{\ast }\tilde{f})\mathrm{lift~}\sigma $ for every $\tilde{f}\in
C^{\infty }(Q)$.  

We assume that the action of $G$ on $(P,\omega )$ preserves $D$. Moreover,
we assume that the momentum map $J:P\rightarrow \mathfrak{g}^{\ast }$ is
constant along $D$. In this case, for each $\xi \in \mathfrak{g}$ and $%
\sigma \in S_{F}^{\infty }(L),$%
\[
\mathcal{Q}_{J_{\xi }}\sigma =J_{\xi }\sigma .
\]%
Hence, a section $\sigma \in S_{F}^{\infty }(L)$ is $G$-invariant if and
only if $J\sigma =0$ and, therefore, the support of $\sigma $ is contained
in $J^{-1}(0)$. We assume here that $J^{-1}(0)$ is a proper closed subset of 
$P$ with empty interior. In this case, only the zero section in $%
S_{F}^{\infty }(L)$ is $G$-invariant, and we are interested in $G$-invariant
generalized sections.

We assume that there is a volume form $dV$ on $Q$ and denote by $\mathcal{%
\tilde{H}}$ the Hilbert space of sections of $\tilde{L}$ which are square
integrable with respect to $dV$. Using the isomorphism $lift:S^{\infty }(%
\tilde{L})\rightarrow S_{F}^{\infty }(L)$, we can lift the scalar product on 
$S^{\infty }(\tilde{L})$ to a scalar product $(\cdot \mid \cdot )$ on $%
S_{F}^{\infty }(L).$ Consider the space 
\[
\mathcal{D}=\{\sigma \in S_{F}^{\infty }(L)\mid (\sigma \mid \sigma )<\infty
\}.
\]%
Let $H$ be the completion of $\mathcal{D}$ with respect to the norm $\sqrt{%
(\sigma \mid \sigma )}$ and let $\mathcal{D}^{\prime }$ be the topological
dual of $\mathcal{D}$. Clearly, $\mathcal{D}\subset \mathcal{H}\subset 
\mathcal{D}^{\prime }$ and $S_{F}^{\infty }(L)\cap \mathcal{H}=\mathcal{D}$. 

We assume that, for each $\xi \in \mathfrak{g}$, the operator $\mathcal{Q}%
_{J_{\xi }}$ on $S_{F}^{\infty }(L)$ preserves $\mathcal{D}$. Hence it
extends to a self adjoint operator on $\mathcal{H}$ and gives rise to a dual
operator $\mathcal{Q}_{J_{\xi }}^{\prime }$ on $\mathcal{D}^{\prime }$ such
that, for every $\xi \in \mathfrak{g}$, $\varphi \in \mathcal{D}^{\prime }$
and $\sigma \in \mathcal{D}$, 
\[
\langle \mathcal{Q}_{J_{\xi }}^{\prime }\varphi \mid \sigma \rangle =\langle
\varphi \mid \mathcal{Q}_{J_{\xi }}\sigma \rangle \text{,}
\]%
where $\langle \mathcal{\cdot }\mid \cdot \rangle $ denotes the evaluation
map. The space of generalized $G$-invariant sections is 
\[
\ker \mathcal{Q}_{J}^{\prime }=\{\varphi \in \mathcal{D}^{\prime }\mid 
\mathcal{Q}_{J_{\xi }}^{\prime }\varphi =0\text{ for all }\xi \in \mathfrak{%
g\}.}
\]%
On the other hand, the range of $\mathcal{Q}_{J}$ in $\mathcal{D}$ is 
\begin{eqnarray}
\mathrm{range~}\mathcal{Q}_{J} &=&\{\dsum\limits_{j=1}^{k}\mathcal{Q}%
_{J_{\xi _{j}}}\sigma _{j}\mid \sigma _{1},...,\sigma _{k}\in \mathcal{D}\}
\label{range} \\
&=&\{\dsum\limits_{j=1}^{k}J_{\xi _{j}}\sigma _{j}\mid \sigma
_{1},...,\sigma _{k}\in \mathcal{D}\}=\mathcal{JD},
\end{eqnarray}%
where $(\xi _{1},...,\xi _{k})$ is a basis in $\mathfrak{g}$.

\begin{lemma}
\label{lemma2}The class $[\sigma ]\in S^{\infty }(L)/\mathcal{J}S^{\infty
}(L)$ of $\sigma \in S^{\infty }(L)$ is uniquely determined by the
restriction of $\sigma $ to any open set containing $J^{-1}(0).$
\end{lemma}

\begin{proof}
Suppose that $(\sup \mathrm{port~}\sigma )\cap J^{-1}(0)=\emptyset $. Hence,
for each $p\in \sup \mathrm{port~}\sigma $, there exists a neighbourhood $%
U_{p}$ of $p$ and $\xi _{p}\in \mathfrak{g}$ such that $J_{\xi _{p}}$ does
not vanish on $U_{p}$. Therefore, we can write%
\[
\sigma _{\mid U_{p}}=J_{\xi _{p}\mid U_{p}}(J_{\xi _{p}\mid
U_{p}})^{-1}\sigma _{\mid U_{p}}=\sum_{j=1}^{k}c_{p}^{j}J_{\xi _{j}\mid
U_{p}}(J_{\xi _{p}\mid U_{p}})^{-1}\sigma _{\mid U_{p}},
\]%
where $c_{p}^{j}$ are components of $\xi _{p}$ with respect to a basis $(\xi
_{1},...,\xi _{k})$ of $\mathfrak{g}$.

The family of open sets $\{U_{p}\mid p\in \sup \mathrm{port~}\sigma \}$
together with the complement $U_{0}$ of $\sup \mathrm{port~}\sigma $ form a
cover of $P$. Let $\{U_{p_{\alpha }}\mid \alpha \in A\}$ be a locally finite
subcover of $P$ and $\{f_{\alpha }\mid \alpha \in A\}$ a corresponding
partition of unity on $P$. Since $\sup \mathrm{port~}f_{\alpha }\subseteq
U_{p_{\alpha }}$, it follows that $f_{\alpha }J_{\xi _{j}\mid U_{p_{\alpha
}}}=J_{\xi _{j}}f_{\alpha }$. Therefore, 
\begin{eqnarray*}
\sigma  &=&\sum_{\alpha \in A}f_{\alpha }\sigma _{\mid U_{p_{\alpha
}}}=\sum_{\alpha \in A}f_{\alpha }\sum_{j=1}^{k}c_{p_{\alpha }}^{j}J_{\xi
_{j}\mid U_{p_{\alpha }}}(J_{\xi _{p_{\alpha }}\mid U_{p_{\alpha
}}})^{-1}\sigma _{\mid U_{p_{\alpha }}} \\
&=&\sum_{j=1}^{k}J_{\xi _{j}}\sum_{\alpha \in A}f_{\alpha }c_{p_{\alpha
}}^{j}(J_{\xi _{p_{\alpha }}\mid U_{p_{\alpha }}})^{-1}\sigma _{\mid
U_{p_{\alpha }}}\in \mathcal{J}S^{\infty }(L)
\end{eqnarray*}%
because $\sum_{\alpha \in A}f_{\alpha }c_{p_{\alpha }}^{j}(J_{\xi
_{p_{\alpha }}\mid U_{p_{\alpha }}})^{-1}\sigma _{\mid U_{p_{\alpha }}}$ are
smooth sections of $L$ for $j=1,...,k$. Hence, $[\sigma ]=0.$

Suppose now $\sigma $ and $\sigma ^{\prime }$ in $S^{\infty }(L)$ are such
that $\sigma $ coincides with $\sigma ^{\prime }$ in an open set $U$
containing $J^{-1}(0)$. Then $(\sigma -\sigma ^{\prime })_{\mid U}=0$ and,
by the argument above, $[\sigma -\sigma ^{\prime }]=0$. Hence, $[\sigma
]=[\sigma ^{\prime }]$.
\end{proof}

We can now state and prove the main result of this subsection.

\begin{theorem}
\label{theorem4}Suppose that $\vartheta (J^{-1}(0))$ has compact closure.
Then, under the assumptions made above,  
\[
\left( S_{F}^{\infty }(L)/\mathcal{J}S^{\infty }(L)\right) ^{G}=\mathcal{D}/%
\mathrm{range~}\mathcal{Q}_{J}\text{.}
\]
\end{theorem}

\begin{proof}
Since $\vartheta :P\rightarrow Q$ is a submersion, every point $p\in P$ is
in the image of a local section $\tau $ of $\vartheta $. Let $\tau _{\alpha
}:U_{\alpha }\rightarrow P$ be a family of sections of $\vartheta $ labelled
by an index set $A$ such that their domains form a locally finite covering
of $Q,$ and let $\{f_{\alpha }\}$ be a partition of unity on $Q$ such that $%
\sup \mathrm{port}~f_{\alpha }\subseteq U_{\alpha }$ for each $\alpha \in A$.

Consider $\sigma \in S_{F}^{\infty }(L)$ such that $[\sigma ]=0$. Then 
\begin{equation}
\sigma =\sum_{j=1}^{k}J_{\xi j}\sigma _{j}  \label{sigma}
\end{equation}%
for some $\sigma _{1},...,\sigma _{k}\in S^{\infty }(L)$. We want to show
that the assumption that $\sigma \in S_{F}^{\infty }(L),$ we may choose
sections $\sigma _{1},...,\sigma _{k}$ in $S_{F}^{\infty }(L)$. Equation (%
\ref{sigma}) implies 
\begin{equation}
\lambda \circ \sigma \circ \tau _{\alpha }=\sum_{j=1}^{k}(J_{\xi _{j}}\circ
\tau _{\alpha })\lambda \circ \sigma _{j}\circ \tau _{\alpha }  \label{alpha}
\end{equation}%
for each $\alpha \in A$. Since $\tau _{\alpha }$ is a section of $\vartheta $%
, it follows that $\lambda \circ \sigma \circ \tau _{\alpha }=\tilde{\sigma}%
\circ \vartheta \circ \tau _{\alpha }=\tilde{\sigma}_{\mid U_{\alpha }}$.
Moreover, since $J$ is constant along fibres of $\vartheta $, it defines a
smooth map $\tilde{J}:Q\rightarrow \mathfrak{g}^{\ast }$ such that $J=\tilde{%
J}\circ \vartheta $. Finally,%
\[
\tilde{\sigma}_{j}^{\prime }=\sum_{\alpha \in A}f_{\alpha }(\lambda \circ
\sigma _{j}\circ \tau _{\alpha })
\]%
is a smooth section of $\tilde{\pi}:\tilde{L}\rightarrow Q$ for each $%
j=1,...,k.$ Hence, for each $q\in Q,$%
\begin{eqnarray*}
\tilde{\sigma}(q) &=&\sum_{\alpha \in A}f_{\alpha }(q)\tilde{\sigma}%
(q)=\sum_{\alpha \in A}f_{\alpha }(q)\sum_{j=1}^{k}J_{\xi _{j}}(\tau
_{\alpha }(q))(\lambda \circ \sigma _{j}\circ \tau _{\alpha })(q) \\
&=&\sum_{j=1}^{k}\tilde{J}_{\xi _{j}}(q)\sum_{\alpha \in A}f_{\alpha
}(q)(\lambda \circ \sigma _{j}\circ \tau _{\alpha })=\sum_{j=1}^{k}\tilde{J}%
_{\xi _{j}}(q)\tilde{\sigma}_{j}^{\prime }(q).
\end{eqnarray*}%
Therefore, 
\begin{equation}
\tilde{\sigma}=\sum_{j=1}^{k}\tilde{J}_{\xi _{j}}\tilde{\sigma}_{j}^{\prime }%
\text{.}  \label{sigma1}
\end{equation}%
For each $j=1,...,k$, the horizontal lift $\sigma _{j}^{\prime }$ of $\tilde{%
\sigma}_{j}^{\prime }$ is in $S_{F}^{\infty }(L)$. By construction, $\sigma =%
\mathrm{lift}~\tilde{\sigma}$.T Hence, $\sigma =\sum_{j=1}^{k}(\tilde{J}%
_{\xi _{j}}\circ \vartheta )\sigma _{j}^{\prime }=\sum_{j=1}^{k}J_{\xi
_{j}}\sigma _{j}^{\prime }$, where $\sigma _{j}^{\prime }\in S_{F}^{\infty
}(L)$. Thus, we have shown that 
\[
S_{F}^{\infty }(L)/\mathcal{J}S^{\infty }(L)=S_{F}^{\infty }(L)/\mathcal{J}%
S_{F}^{\infty }(L).
\]

Since $\vartheta (J^{-1}(0))$ is compact, there exists an open set $V$ in $Q$
with compact closure $\bar{V}$ such that $\vartheta (J^{-1}(0))\subseteq V$.
For every smooth section $\tilde{\sigma}$ of $\tilde{\pi}:\tilde{L}%
\rightarrow Q$, there exists a compactly supported section $\tilde{\sigma}%
^{\prime }$ which coincides with $\sigma $ on $\vartheta ^{-1}(V)$. The
horizontal lifts $\sigma $ and $\sigma ^{\prime }$ of $\tilde{\sigma}$ and $%
\tilde{\sigma}^{\prime }$, respectively, coincide on $\vartheta ^{-1}(V)$
and Lemma \cite{lemma2} implies that $[\sigma ]=[\sigma ^{\prime }]$. Since $%
\tilde{\sigma}^{\prime }$ is compactly supported, it follows that it is
square integrable over $Q$. Hence, $\sigma ^{\prime }\in S_{F}^{\infty
}(L)\cap \mathcal{H=D}$. Therefore, 
\[
S_{F}^{\infty }(L)/\mathcal{J}S_{F}^{\infty }(L)=\mathcal{D}/\mathcal{J}%
S_{F}^{\infty }(L)\text{,}
\]%
where $\mathcal{D}/\mathcal{J}S^{\infty }(L)$ is the image of $\mathcal{D}$
under the projection $S_{F}^{\infty }(L)\rightarrow S_{F}^{\infty }(L)/%
\mathcal{J}S_{F}^{\infty }(L).$

Suppose now that $\tilde{\sigma}$ in equation (\ref{sigma1}) is compactly
supported. Multiplying equation (\ref{sigma1}) by a compactly supported
function $\tilde{f}$, such that $\tilde{f}_{\mid \sup \mathrm{port}\text{~}%
\tilde{\sigma}}=1,$ we see that equation (\ref{sigma1}) holds with compactly
supported sections $\tilde{f}\tilde{\sigma}_{1}^{\prime },...,\tilde{f}%
\tilde{\sigma}_{k}^{\prime }$ on the right hand side. This implies that 
\[
\mathcal{D}/\mathcal{J}S_{F}^{\infty }(L)=\mathcal{D}/\mathcal{JD}\text{.}
\]

Finally, consider  the space $(S_{F}^{\infty }(L)/\mathcal{J}S_{F}^{\infty
}(L))^{G}$ of $G$-invariant elements in $S_{F}^{\infty }(L)/\mathcal{J}%
S_{F}^{\infty }(L)$. By definition, it consists of $[\sigma ]\in
S_{F}^{\infty }(L)/\mathcal{J}S_{F}^{\infty }(L)$ such that $\mathcal{Q}%
_{J_{\xi }}\sigma \in \mathcal{J}S_{F}^{\infty }(L)$ for each $\xi \in 
\mathfrak{g}$. However, for our polarization $F$, each operator $\mathcal{Q}%
_{J_{\xi }}$ is an operator of multiplication by $J_{\xi }$. In other words, 
$\mathcal{Q}_{J_{\xi }}\sigma =J_{\xi }\sigma $ for each $\sigma \in
S_{F}^{\infty }(L)$. Taking into account equation (\ref{range}), we obtain 
\[
(S_{F}^{\infty }(L)/\mathcal{J}S^{\infty }(L))^{G}=S_{F}^{\infty }(L)/%
\mathcal{J}S_{F}^{\infty }(L)=\mathcal{D}/\mathcal{JD}=\mathcal{D}/\mathrm{%
range~}\mathcal{Q}_{J}
\]%
which completes the proof.
\end{proof}

\begin{corollary}
\label{corollary1}There is a natural duality between $(S_{F}^{\infty }(L)/%
\mathcal{J}S^{\infty }(L))^{G}$ and $\ker \mathcal{Q}_{J}^{\prime }\subset 
\mathcal{D}^{\prime }.$
\end{corollary}

\begin{proof}
There is a natural duality between $\ker \mathcal{Q}_{J}^{\prime }\subset 
\mathcal{D}^{\prime }$ and $\mathcal{D}/\mathrm{range~}\mathcal{Q}_{J}$ such
that, the evaluation of $\varphi \in \ker \mathcal{Q}_{J}^{\prime }$ on the
class $\left\{ \sigma \right\} \in \mathcal{D}/\mathrm{range~}\mathcal{Q}_{J}
$ of $\sigma \in \mathcal{D}$ is given by%
\[
\langle \varphi \mid \left\{ \sigma \right\} \rangle =\langle \varphi \mid
\sigma \rangle .
\]%
But $(S_{F}^{\infty }(L)/\mathcal{J}S^{\infty }(L))^{G}=\mathcal{D}/\mathrm{%
range}~\mathcal{Q}_{J}$ by Theorem \ref{theorem4}.
\end{proof}

\subsection{Example: the ground state of a free particle.}

In this subsection we give an example of quantization of algebraic reduction
for an improper action of a Lie group.

\subsubsection{Quantization representation}

Consider $P=\mathbb{R}^{2}$, with coordinates $(p,q)$, and a symplectic form 
$\omega =dp\wedge dq.$ The action of $G=\mathbb{R}$ is $\Phi
_{t}(p,q)=(p,q+tp)$ for all $t\in G$ and $(p,q)\in P$. This action is
Hamiltonian with the momentum map $J:P\mapsto \mathfrak{g}=\mathbb{R}$ given
by the Hamiltonian of a free particle with unit mass. In other words, 
\begin{equation}
J(p,q)=\frac{1}{2}p^{2}.  \label{energy}
\end{equation}%
The fixed point set of this action is the zero level set 
\[
J^{-1}(0)=\{(p,q)\mid p=0\}\cong \mathbb{R}.
\]%
Since $\mathbb{R}$ is not compact it follows that the action $\Phi $ is not
proper.

We take $L=P\times \mathbb{C}$ with a trivializing section $\sigma _{0}$
such that 
\[
\sigma _{0}(p,q)=((p,q),1) 
\]%
for all $(p,q)\in P$. The usual Schr\"{o}dinger (position) representation is
obtained by choosing a connection $\nabla $ such that $\nabla \sigma
_{0}=-i\hbar ^{-1}pdq\otimes \sigma _{0},$ and a vertical polarization
spanned by $\frac{\partial }{\partial p}$. However, the vertical
polarization is not preserved by our group action. Therefore, we choose a
momentum representation given by a polarization $F=\mathrm{span}\{\frac{%
\partial }{\partial q}\}\otimes \mathbb{C}$, see \cite{GQ-QM}. Consider a
trivializing section $\sigma _{1}=\exp (i\hbar ^{-1}pq)\sigma _{0}$. We have 
\begin{eqnarray*}
\nabla \sigma _{1} &=&\nabla \{\exp (i\hbar ^{-1}pq)\sigma _{0}\}=d\{\exp
(i\hbar ^{-1}pq)\}\otimes \sigma _{0}+\exp (i\hbar ^{-1}pq)\nabla \sigma _{0}
\\
&=&i\hbar ^{-1}qdp\otimes \exp (i\hbar ^{-1}pq)\sigma _{0}=i\hbar
^{-1}qdp\otimes \sigma _{1}.
\end{eqnarray*}%
In particular, $\nabla _{\frac{\partial }{\partial q}}\sigma _{1}=0$ which
implies that the space of sections of $L$ that are covariantly constant
along $F$ is 
\[
S_{F}^{\infty }(L)=\{\psi (p)\sigma _{1}\mid \psi \in C^{\infty }(\mathbb{R}%
)\}. 
\]

The Poisson algebra of quantizable functions is%
\[
{C_{F}^{\infty }(P)=\{h}_{1}{(p)+qh}_{2}(p)\mid h_{1},h_{2}\in C^{\infty }(%
\mathbb{R})\}. 
\]%
For each ${f}\in C_{F}^{\infty }(P)$, the corresponding quantum operator on $%
S_{F}^{\infty }(L)$ is given by 
\[
\mathcal{Q}_{f}(\psi \sigma _{1})=(-i\hbar \nabla _{X_{f}}+f)(\psi \sigma
_{1})=\{-i\hbar (X_{f}\psi )+(f+\langle qdp\mid X_{f}\rangle )\}\sigma _{1}. 
\]%
Since $X_{h_{1}(p)}=\frac{\partial h_{1}}{\partial p}\frac{\partial }{%
\partial q},$ it follows that $\mathcal{Q}_{h_{1}(p)}$ is the operator of
multiplication by $h_{1}(p)$. In other words, 
\begin{equation}
\mathcal{Q}_{h_{1}(p)}(\psi \sigma _{1})=h_{1}(p)\psi (p)\sigma _{1}.
\label{h1}
\end{equation}%
Similarly, $X_{qh_{2}(p)}=q\frac{\partial h_{2}}{\partial p}\frac{\partial }{%
\partial q}-h_{2}(p)\frac{\partial }{\partial p}$ and 
\begin{equation}
\mathcal{Q}_{qh_{2}(p)}(\psi \sigma _{1})=-i\hbar h_{2}(p)\frac{d\psi (p)}{dp%
}\sigma _{1}.  \label{qh2}
\end{equation}

The momentum map $J=\frac{1}{2}p^{2}$ gives rise to an operator $\mathcal{Q}%
_{J}$ of multiplication by $J$. In other words, 
\begin{equation}
\mathcal{Q}_{J}\sigma =J\sigma  \label{QJ}
\end{equation}%
for all $\sigma \in S_{F}^{\infty }(L)$. For each $t\in \mathbb{R}$,
considered as the Lie algebra of $\mathbb{R}$, the corresponding momentum is 
$J_{t}=tJ$, and the infinitesimal representation $t\mapsto (-i\hbar )^{-1}%
\mathcal{Q}_{tJ}$ integrates to a representation $\mathcal{R}$ of $\mathbb{R}
$ on $S_{F}^{\infty }(L)$ given by%
\begin{equation}
\mathcal{R}_{t}\psi (p)\sigma _{1}=\exp (-i\hbar ^{-1}tJ(p))\psi (p)\sigma
_{1}.  \label{Rt}
\end{equation}

The representation (\ref{Rt}) gives rise to a unitary representation $%
\mathcal{U}$ of $\mathbb{R}$ on 
\[
\mathcal{H}=\{\psi (p)\sigma _{1}\mid \psi \in L^{2}(\mathbb{R})\}
\]%
with a scalar product 
\[
(\psi _{1}(p)\sigma _{1},\psi _{2}(p)\sigma _{1})=\int_{\mathbb{R}}\bar{\psi}%
_{1}(p)\psi _{2}(p)dp.
\]%
For each $t\in \mathbb{R}$, and $\psi (p)\sigma _{1}\in \mathcal{H},$ 
\[
\mathcal{U}_{t}\psi (p)\sigma _{1}=\exp (-i\hbar ^{-1}tJ(p))\psi (p)\sigma
_{1}.
\]

Let $\mathcal{D}=S_{F}^{\infty }(L)\cap \mathcal{H}$ and $\mathcal{D}%
^{\prime }$ be the space of continuous linear functionals on $\mathcal{D}$.
For each $\varphi \sigma _{1}\in \mathcal{D}^{\prime }$ and $\psi \sigma
_{1}\in \mathcal{D}$ we write the evaluation of $\varphi \sigma _{1}$ on $%
\psi \sigma _{1}$ as if it were given by integration, namely%
\[
\langle \varphi \sigma _{1}\mid \psi \sigma _{1}\rangle =\int_{\mathbb{R}%
}\varphi (p)\psi (p)dp.
\]%
The multiplication operator $\psi \sigma _{1}\mapsto J\psi \sigma _{1}$ on $%
\mathcal{H}$ has continuous spectrum $[0,\infty )$. For each eigenvalue $%
\lambda \geq 0$, eigenvectors corresponding to $\lambda $ are generalized
sections $\varphi \sigma _{1}$ satisfying the equation 
\begin{equation}
(p^{2}-2\lambda )\varphi (p)\sigma _{1}=0\text{. }  \label{eigenvectors}
\end{equation}%
If $\lambda >0$, the eigenspace corresponding to $\lambda $ is given by 
\[
\varphi (p)\sigma _{1}=A\delta (p-\sqrt{2\lambda })\sigma _{1}+B\delta (p-%
\sqrt{2\lambda })\sigma _{1},
\]%
where $A$ and $B$ are arbitrary constants and $\delta $ is the Dirac delta
function. The eigenspace corresponding to $\lambda =0$ is the kernel of $%
\mathcal{Q}_{J}^{\prime }$. It is given by 
\begin{equation}
\ker \mathcal{Q}_{J}^{\prime }=\{(A\delta (p)+B\delta ^{\prime }(p))\sigma
_{1}\mid A,B\in \mathbb{C}\}.  \label{kerQJ}
\end{equation}

The ring of smooth $G$-invariant functions on $\mathbb{R}^{2}$ is denoted by 
$C^{\infty }(\mathbb{R}^{2})^{G}$. A function $f(p,q)$ belongs to $C^{\infty
}(\mathbb{R}^{2})^{G}$ if $f(p,q+pt)=f(q,p)$ for all $t$. In other words, $%
f(p,q)$ is in $C^{\infty }(\mathbb{R}^{2})^{G}$ if and only if it is
independent of $q$. Thus, elements $f\in C^{\infty }(\mathbb{R}^{2})^{G}$
are of the form $f=k(p)$ for $k\in C^{\infty }(\mathbb{R}).$ For every $%
f=k(p)\in C^{\infty }(\mathbb{R}^{2})^{G}$, the action of the quantum
operator $\mathcal{Q}_{f}^{\prime }$ on $\mathcal{D}^{\prime }$ preserves $%
\ker \mathcal{Q}_{J}^{\prime }$. For each $(A\delta (p)+B\delta ^{\prime
}(p))\sigma _{1}\in \ker \mathcal{Q}_{J}^{\prime }$, 
\begin{eqnarray}
\mathcal{Q}_{k(p)}^{\prime }(A\delta (p)+B\delta ^{\prime }(p))\sigma _{1}
&=&(Ak\left( p)\delta (p)+Bk(p)\delta ^{\prime }(p)\right) \sigma _{1}
\label{Qf(p)} \\
&=&(Ak\left( 0)-Bk^{\prime }(0))\delta (p)+Bk(0)\delta ^{\prime }(p)\right)
\sigma _{1}
\end{eqnarray}

\subsubsection{Reduction at 0}

In our example, 
\[
\mathcal{J}=\{f\in C^{\infty }(\mathbb{R}^{2})\mid f=p^{2}\tilde{f}\text{
for }\tilde{f}\in C^{\infty }(\mathbb{R}^{2})\}. 
\]%
For $f\in C^{\infty }(\mathbb{R}^{2})$, we denote by $[f]$ the equivalence
class of $f$ in $C^{\infty }(\mathbb{R}^{2})/\mathcal{J}$. Using the Taylor
expansion of $f(p,q)$ with respect to $p$ we can see that $[f]$ is uniquely
described by the first two terms $f(0,q)+f_{p}(0,q)p$, where $f_{p}=\frac{%
\partial f}{\partial p}.$ Hence, we can make an identification 
\begin{eqnarray*}
C^{\infty }(\mathbb{R}^{2})/\mathcal{J} &\equiv &\{f(0,q)+f_{p}(0,q)p\mid
f\in C^{\infty }(\mathbb{R}^{2})\} \\
&\equiv &\{h_{1}(q)+ph_{2}(q)\in C^{\infty }(\mathbb{R}^{2})\mid
h_{1},h_{2}\in C^{\infty }(\mathbb{R})\}.
\end{eqnarray*}%
Under this identification, $C^{\infty }(\mathbb{R}^{2})/\mathcal{J}$ is
interpreted as a subspace of $C^{\infty }(\mathbb{R}^{2})$ and the
projection map $C^{\infty }(\mathbb{R}^{2})\rightarrow C^{\infty }(\mathbb{R}%
^{2})/\mathcal{J}:f\mapsto \lbrack f]$ is a retraction.

The action of $\mathbb{R}$ on $C^{\infty }(\mathbb{R}^{2})/\mathcal{J}$ is
given by $\tilde{\Phi}_{t}^{\ast }[f]=[\Phi _{t}^{\ast }f]$ for every $t\in 
\mathbb{R}$, and $[f]\in C^{\infty }(\mathbb{R}^{2})/\mathcal{J}$. Since $%
(\Phi _{t}^{\ast }f)(p,q)=f(p,q+tp)$, it follows that 
\[
\tilde{\Phi}_{t}^{\ast }(h_{1}(q)+ph_{2}(q))=h_{1}(q+pt)+ph_{2}(q+tp). 
\]%
Differentiating with respect to $t$ and setting $t=0$, we get 
\[
\frac{d}{dt}\tilde{\Phi}_{t}^{\ast }(h_{1}(q)+ph_{2}(q))_{\mid t=0}=p\frac{%
dh_{1}(q)}{dq}+p^{2}\frac{dh_{2}(q)}{dq}. 
\]%
Thus, an element $h_{1}(q)+ph_{2}(q)$ of $C^{\infty }(\mathbb{R}^{2})/%
\mathcal{J}$ is invariant if $h_{1}$ is constant. Hence, the reduced algebra
can be presented as follows. 
\[
(C^{\infty }(\mathbb{R}^{2})/\mathcal{J})^{G}\equiv \{a+ph(q)\mid a\in 
\mathbb{R},\text{ }h\in C^{\infty }(\mathbb{R})\}\subset C^{\infty }(\mathbb{%
R}^{2}). 
\]%
On the other hand, 
\[
C^{\infty }(\mathbb{R}^{2})^{G}/\mathcal{J}\equiv \{a+bp\mid a,b\in \mathbb{R%
}\}\subset C^{\infty }(\mathbb{R}^{2}). 
\]

For the polarization $F=\mathrm{span}\{\frac{\partial }{\partial q}\}\otimes 
\mathbb{C}$, the Poisson algebra of quantizable functions is ${C_{F}^{\infty
}(}\mathbb{R}^{2}{)=\{f}_{1}{(p)+qf}_{2}(p)\mid f_{1},f_{2}\in C^{\infty }(%
\mathbb{R})\}$. Hence, 
\begin{eqnarray*}
(C_{F}^{\infty }(\mathbb{R}^{2})/\mathcal{J})^{G} &=&\{[f]\in (C^{\infty }(%
\mathbb{R}^{2})/\mathcal{J})^{G}\mid f\in C_{F}^{\infty }(P)\} \\
&\equiv &\{a+bp+cpq\mid a,b,c\in \mathbb{R}\}\subset C^{\infty }(\mathbb{R}%
^{2}).
\end{eqnarray*}

Next, we identify the modules involved in quantization. We have 
\[
\mathcal{J}S^{\infty }(L)=\{p^{2}\psi (q,p)\sigma _{1}\mid \psi \in
C^{\infty }(\mathbb{R}^{2})\otimes \mathbb{C}\}. 
\]%
As before, we can identify the quotient $S^{\infty }(L)/\mathcal{J}S^{\infty
}(L)$ as a submodule of $S^{\infty }(L)$ such that the projection map $%
S^{\infty }(L)\rightarrow S^{\infty }(L)/\mathcal{J}S^{\infty }(L):\sigma
\mapsto \lbrack \sigma ]$ is a retraction. This identification is given by 
\[
S^{\infty }(P,L)/\mathcal{J}S^{\infty }(P,L)\equiv \{(\chi _{1}(q)+p\chi
_{2}(q))\sigma _{1}\mid \chi _{1},\chi _{2}\in C^{\infty }(\mathbb{R}%
)\otimes \mathbb{C}\}\subset S^{\infty }(L). 
\]%
Similarly, 
\[
S_{F}^{\infty }(L)=\{\psi (p)\sigma _{1}\mid \psi \in C^{\infty }(\mathbb{R}%
)\otimes \mathbb{C\}}\subset S^{\infty }(L), 
\]%
and we can make an identification 
\[
S_{F}^{\infty }(L)/\mathcal{J}S^{\infty }(L)\equiv \{(d+ep)\sigma _{1}\mid
d,e\in \mathbb{C}\}\subset S_{F}^{\infty }(L)\}. 
\]%
It is easy to see that elements of $S_{F}^{\infty }(L)/\mathcal{J}S^{\infty
}(L)$ are $G$-invariant. In other words, 
\[
S_{F}^{\infty }(L)/\mathcal{J}S^{\infty }(L)=(S_{F}^{\infty }(L)/\mathcal{J}%
S^{\infty }(L))^{G}. 
\]%
If $[f]=a+bp+cpq$ and $[\sigma ]=(d+ep)\sigma _{1}$, then equations (\ref{h1}%
), (\ref{qh2}), and $\mathcal{Q}_{[f]}[\sigma ]=[\mathcal{Q}_{f}\sigma ]$
give 
\[
\mathcal{Q}_{[a+bp+cpq]}[(d+ep)\sigma _{1}]=\{ad+(bd+ae-i\hbar ce)p\}\sigma
_{1}. 
\]%
In particular, if $f=k(p)\in C^{\infty }(\mathbb{R}^{2})^{G}$ then $%
[k(p)]=k(0)+k^{\prime }(0)p$, and 
\begin{equation}
\mathcal{Q}_{[k(p)]}[(d+ep)\sigma _{1}]=\{k(0)d+(k^{\prime
}(0)d+k(0)e)p\}\sigma _{1}.  \label{Q[f(p)]}
\end{equation}

By Corollary \ref{corollary1} there is a duality between $(S_{F}^{\infty
}(L)/\mathcal{J}S^{\infty }(L))^{G}=\{(d+ep)\sigma _{1}\mid d,e\in \mathbb{C}%
\}$ and $\ker \mathcal{Q}_{J}^{\prime }=\{(A\delta (p)+B\delta ^{\prime
}(p))\sigma _{1}\mid A,B\in \mathbb{C}\}$. It is given by 
\[
\langle (A\delta (p)+B\delta ^{\prime }(p))\sigma _{1}\mid (d+ep)\sigma
_{1}\rangle =Ad-Be. 
\]

\section{Non-zero co-adjoint orbits}

\subsection{Quantization}

We turn now to quantization of algebraic reduction at a non-zero coadjoint
orbit $O.$ Recall that algebraic reduction at $O$ is defined as algebraic
reduction at $0\in \mathfrak{g}^{\ast }$ of the action of $G$ on $(P\times
O,pr_{1}^{\ast }\omega -pr_{2}^{\ast }\omega _{O})$. In other words, the
Poisson algebra of algebraic reduction at $O$ is $(C^{\infty }(P\times O)/%
\mathcal{J}_{P\times O})^{G}$, where $\mathcal{J}_{P\times O}$\ is the ideal
in $C^{\infty }(P\times O)$\ generated by components of the momentum map $%
J_{P\times O}:P\times O\rightarrow \mathfrak{g}^{\ast }$ given by (\ref{JPO}%
). Hence, quantization of algebraic reduction at $O$ is given by
quantization of algebraic reduction of $(P\times O,pr_{1}^{\ast }\omega
-pr_{2}^{\ast }\omega _{O})$ at the zero value of $J_{P\times O}$.

We assume that $(O,\omega _{O})$ is a quantizable coadjoint orbit. Let $\pi
_{O}:L_{O}\rightarrow O$ be a prequantization line bundle for $(O,\omega
_{O})$ and $\nabla _{O}$ a connection on $L_{O}$ with curvature $\omega _{O}$%
. We assume that there is a $G$-invariant polarization $F_{O}$ of $(O,\omega
_{O})$. To simplify the notation, we denote by $C_{F}^{\infty }(O)$ the
space of functions in $C^{\infty }(O)$ such that their Hamiltonian vector
fields preserve the polarization $F_{O}$. Similarly, we denote by $%
S_{F}^{\infty }(L_{O})$ the space of smooth sections of $L_{O}$ which are
covariantly constant along $F_{O}$. Quantization associates to each $%
f_{O}\in C_{F}^{\infty }(O)$ a linear operator $(\mathcal{Q}_{O})_{f_{O}}$
on the space $S_{F}^{\infty }(L_{O})$ given by 
\begin{equation}
(\mathcal{Q}_{O})_{f_{O}}\sigma _{O}=(-i\hbar (\nabla
_{O})_{X_{f_{O}}}+f_{O})\sigma _{O}  \label{Quantization1}
\end{equation}%
for every $\sigma _{O}\in S_{F}^{\infty }(O)$. Here, $X_{f_{O}}$ denotes the
Hamiltonian vector field on $O$ defined in terms of the symplectic form $%
\omega _{O}$. We denote by $\mathcal{R}_{O}$ the representation of $G$ on $%
S_{F}^{\infty }(O)$ obtained by integration of the representation $\xi
\mapsto (-i\hbar )^{-1}(\mathcal{Q}_{O})_{J_{\xi }}$ of $\mathfrak{g}$.

Complex conjugation $z\mapsto \bar{z}$ in $L_{O}$ is an automorphism of a
real vector bundle over $O$, but it conjugates the complex structure. We
denote by $\bar{L}_{O}$ the complex line bundle over $O$ with the conjugate
complex structure. If $\sigma _{O}$ is a non-zero section of $L_{O}$, then $%
\bar{\sigma}_{O}$ is a non-zero section of $\bar{L}_{O}$. Moreover, $\nabla
_{O}\sigma _{O}=-i\hbar ^{-1}\theta _{O}\otimes \sigma _{O}$ implies $\bar{%
\nabla}_{O}\bar{\sigma}_{O}=i\hbar ^{-1}\bar{\theta}_{O}\otimes \bar{\sigma}%
_{O}$. Thus, the pull back by $\bar{\sigma}_{O}$ of the connection form of $%
\bar{\nabla}_{O}$ is given by $\bar{\theta}_{O}=-\theta _{O}$. Hence, $%
-\omega _{O}$ is the curvature of the $\bar{\nabla}_{O}$ on $\bar{L}_{O}$.
Moreover, for $f_{O}\in C^{\infty }(O)$, the Hamiltonian vector field of $%
f_{O}$ relative to the symplectic form $-\omega _{O}$ is $\bar{X}%
_{f_{O}}=-X_{f_{O}}$, where $X_{f_{O}}$ is the Hamiltonian vector field with
respect to $\omega _{O}$. Therefore, for every $\bar{\sigma}_{O}\in
S^{\infty }(\bar{L}_{O})$ 
\[
(\bar{\nabla}_{O})_{\bar{X}_{f_{O}}}\bar{\sigma}_{O}=i\hbar ^{-1}\langle 
\bar{\theta}_{O}\mid \bar{X}_{f_{O}}\rangle \bar{\sigma}_{O}=(\overline{%
-i\hbar ^{-1}\langle \theta _{O}\mid X_{f_{O}}\rangle \sigma })=\overline{%
(\nabla _{O})_{X_{f_{O}}}\sigma _{O}}. 
\]

By assumption, $F_{O}$ is a positive polarization of $(O,\omega _{O})$, that
is, $i\omega _{O}(w,\bar{w})\geq 0$ for all $w\in F_{O}$. This implies that $%
\bar{F}_{O}$ is a positive polarization of $(O,-\omega _{O})$. Thus, the
chosen quantization structure on $(O,\omega _{O})$ induces a quantization
structure on $(O,-\omega _{O})$ given by $\bar{L}_{O},$ $\bar{\nabla}_{O}$, $%
\bar{F}_{O}$. We denote by $S_{\bar{F}}^{\infty }(\bar{L}_{O})$ the space of
smooth sections of $\bar{L}$ which are covariantly constant along $\bar{F}%
_{O}$. Similarly, we denote by $C_{\bar{F}}^{\infty }(O)$ the Poisson
algebra of functions in $C^{\infty }(O)$ such that their Hamiltonian vector
fields with respect to the symplectic form $(-\omega _{O})$ preserve the
polarization $\bar{F}_{O}$. Since functions in $C^{\infty }(O)$ are real it
follows that $C_{\bar{F}}^{\infty }(O)=C_{F}^{\infty }(O)$. 

Let $\overline{\mathcal{Q}}_{O}:C_{\bar{F}}^{\infty }(O)\times S_{\bar{F}%
}^{\infty }(\bar{L}_{O})\rightarrow S_{\bar{F}}^{\infty }(\bar{L}_{O})$ be
the quantization map corresponding to the quantiation structure $\bar{L}_{O},
$ $\bar{\nabla}_{O}$, $\bar{F}_{O}$. For every $\bar{\sigma}_{O}\in S_{\bar{F%
}}^{\infty }(\bar{L}_{O})$, its complex conjugate $\sigma _{O}=\overline{%
\bar{\sigma}_{O}}$ is in $S_{F}^{\infty }(L_{O})$. Hence, 
\[
\overline{(\mathcal{Q}_{O})_{f_{O}}\sigma _{O}}=\overline{((-i\hbar (\nabla
_{O})_{X_{f_{O}}}+f_{O})\sigma _{O})}=(-i\hbar (\bar{\nabla}_{O})_{\bar{X}%
_{f_{O}}}+f_{O})\bar{\sigma}_{O}=(\overline{\mathcal{Q}}_{O})_{f_{O}}\bar{%
\sigma}_{O}.
\]%
We denote by $\overline{\mathcal{R}}_{O}$ the representation of $G$ on $S_{%
\bar{F}}^{\infty }(O)$ obtained by integration of the representation $\xi
\mapsto (-i\hbar )^{-1}\overline{\mathcal{Q}}_{J_{\xi }}$ of $\mathfrak{g}$.

We consider a line bundle $\pi _{P\times O}:L_{P\times O}\rightarrow P\times
O$ defined as the tensor product of the bundles $\pi :L\rightarrow P$ and $%
\bar{\pi}_{O}:\bar{L}_{O}\rightarrow O$. More precisely, 
\[
L_{P\times O}=pr_{1}^{\ast }L\otimes pr_{2}^{\ast }\bar{L}_{O},
\]%
where $pr_{1}^{\ast }L$ and $pr_{2}^{\ast }\bar{L}_{O}$ are the pullbacks of 
$L$ and $\bar{L}_{O}$ by the projection maps on the first and the second
factor of $P\times O$, respectively. Local sections of $\pi _{P\times
O}:L_{P\times O}\rightarrow P\times O$ are linear combinations of sections
of the form $\sigma _{P\times O}=\sigma \otimes \bar{\sigma}_{O}$, where $%
\sigma $ and $\bar{\sigma}_{O}$ are local sections of $L$ and $\bar{L}_{O}$,
respectively. Let $\nabla _{P\times O}$ be a connection on $L_{P\times O}$
defined by 
\[
\nabla _{P\times O}(\sigma \otimes \bar{\sigma}_{O})=\nabla \sigma \otimes 
\bar{\sigma}_{O}+\sigma \otimes \bar{\nabla}_{O}\bar{\sigma}_{O}.
\]%
The connection $\nabla _{P\times O}$ satisfies the prequantization condition
for $\omega _{P\times O}=pr_{1}^{\ast }\omega -pr_{2}^{\ast }\omega _{O}$.

We choose the polarization $F_{P\times O}$ to be the direct sum $F\oplus 
\bar{F}_{O}$. It is a strongly admissible positive $G$-invariant
polarization of $(P\times O,\omega _{P\times O})$. The space of smooth
functions on $P\times O$ such that their Hamiltonian vector fields preserve
the polarization $F_{P\times O}=F\oplus \bar{F}_{O}$ will be denoted by $%
C_{F}^{\infty }(P\times O)$. A section $\sigma \otimes \bar{\sigma}_{O}$ of $%
L_{P\times O}$ is covariantly constant along $F_{P\times O}$ if and only if $%
\sigma $ is covariantly constant along $F$ and $\bar{\sigma}_{O}$ is
covariantly constant along $\bar{F}_{O}$. Therefore, the space of smooth
sections of $L_{P\times O}$ which are covariantly constant along $F_{P\times
O}$ is 
\[
S_{F}^{\infty }(L_{P\times O})=S_{F}^{\infty }(L)\otimes S_{\bar{F}}^{\infty
}(\bar{L}_{O}).
\]%
For each $f\in C_{F}^{\infty }(P)$, and $f_{O}\in C_{F}^{\infty }(O)$, the
quantization operators $(\mathcal{Q}_{P\times O})_{pr_{1}^{\ast }f}$ and $(%
\mathcal{Q}_{P\times O})_{pr_{2}^{\ast }f_{O}}$ on $S^{\infty }(L_{P\times
O})$ are given by 
\begin{eqnarray*}
(\mathcal{Q}_{P\times O})_{pr_{1}^{\ast }f}(\sigma \otimes \bar{\sigma}_{O})
&=&(\mathcal{Q}_{f}\sigma )\otimes \bar{\sigma}_{O}, \\
(\mathcal{Q}_{P\times O})_{pr_{2}^{\ast }f_{O}}(\sigma \otimes \bar{\sigma}%
_{O}) &=&\sigma \otimes (\overline{\mathcal{Q}}_{O})_{f_{O}}\bar{\sigma}_{O}.
\end{eqnarray*}%
Thus, the quantization representation $\mathcal{R}_{P\times O}$ of $G$ on $%
S_{F}^{\infty }(L_{P\times O})$ is the tensor product of the quantization
representations $\mathcal{R}$ $\,$and $\overline{\mathcal{R}}_{O}$.

Quantization of algebraic reduction of $(P\times O,\omega _{P\times O})$ at $%
0\in \mathfrak{g}^{\ast }$ is interpreted as the quantization of algebraic
reduction of $(P,\omega )$ at $O$. Hence, the Poisson algebra of quantizable
elements of $(C^{\infty }(P\times O)/\mathcal{J}_{P\times O})^{G}$ in the
polarization $F_{P\times O}$ is $(C_{F}^{\infty }(P\times O)/\mathcal{J}%
_{P\times O})^{G}$. The representation space of quantization of algebraic
reduction at $O$ is $\left( S_{F}^{\infty }(L_{P\times O})/(\mathcal{J}%
_{P\times O}S^{\infty }(L_{P\times O}))\right) ^{G}$. Finally, quantization
assigns to each $[f_{P\times O}]\in (C_{F}^{\infty }(P\times O)/\mathcal{J}%
_{P\times O})^{G}$ an operator $\mathcal{Q}_{[f_{P\times O}]}$ on $%
(S_{F}^{\infty }(L_{P\times O})/\mathcal{J}_{P\times O}S^{\infty
}(L_{P\times O}))^{G}$ such that%
\[
\mathcal{Q}_{[f_{P\times O}]}[\sigma _{P\times O}]=[(\mathcal{Q}_{P\times
O})_{f_{P\times O}}\sigma _{P\times O}]. 
\]

\subsection{K\"{a}hler polarizations revisited}

We have seen that, for a K\"{a}hler polarization $F$, quantization of
algebraic reduction at $0\in \mathfrak{g}^{\ast }$ provides complete
information about the space of $G$-invariant polarized sections of $L$. In
this subsection we are going to investigate to what extent quantization of
algebraic reduction at a non-zero quantizable coadjoint orbit $O$ provides
information about sections on which the quantization representation is
equivalent to the irreducible representation corresponding to $O$.

We assume that the polarizations $F$ and $F_{O}$ are K\"{a}hler, which
implies that $F_{P\times O}=F\oplus \bar{F}_{O}$ is a K\"{a}hler
polarization of $(P\times O,pr_{1}^{\ast }\omega -pr_{2}^{\ast }\omega _{O})$%
. In addition, we assume that quantizations of $(P,\omega ),$ $(O,\omega
_{O}),$ $(O,-\omega _{O})$ and $(P\times O,pr_{1}^{\ast }\omega
-pr_{2}^{\ast }\omega _{O})$ in terms of the chosen polarizations give rise
to unitary representations $\mathcal{U}$, $\mathcal{U}_{O}$, $\overline{%
\mathcal{U}}_{O}$ and $\mathcal{U}_{P\times O}$ of $G$ on Hilbert spaces $%
\mathcal{H}$, $\mathcal{H}_{O},$ $\overline{\mathcal{H}}_{O}$ and $\mathcal{H%
}_{P\times O}=\mathcal{H}\otimes \overline{\mathcal{H}}_{O}$, respectively.
Here, $\overline{\mathcal{H}}_{O}$ is the Hilbert space of square integrable
sections of $\bar{L}_{O}$ which are covariantly constant along $\bar{F}_{O}$%
. This implies that the representation $\mathcal{U}_{P\times O}$ of $G$ on $%
\mathcal{H}_{P\times O}$ is the tensor product of the representation $%
\mathcal{U}$ on $\mathcal{H}$ and the representation $\overline{\mathcal{U}}%
_{O}$ on $\overline{\mathcal{H}}_{O}$. We assume also that the
representation $\mathcal{U}_{O}$ is irreducible.

There is a natural bilinear pairing of $\mathcal{H}_{O}$ and $\overline{%
\mathcal{H}}_{O}$ given as follows. For $\sigma _{O}\in \mathcal{H}_{O}$ and 
$\bar{\sigma}_{O}^{\prime }\in \overline{\mathcal{H}}_{O}$, the complex
conjugate of $\bar{\sigma}_{O}^{\prime }$ is $\sigma _{O}^{\prime }\in 
\mathcal{H}_{O}.$ We define the pairing $\langle \sigma _{O}\mid \bar{\sigma}%
_{O}^{\prime }\rangle $ to be the scalar product $(\sigma _{O}\mid \sigma
_{O}^{\prime })_{\mathcal{H}_{O}}$ of $\sigma _{O}$ and $\sigma _{O}^{\prime
}$ in $\mathcal{H}_{O}.$ In other words, 
\[
\langle \cdot \mid \cdot \rangle :\mathcal{H}_{O}\times \overline{\mathcal{H}%
}_{O}\rightarrow \mathbb{C}:(\sigma _{O},\bar{\sigma}_{O}^{\prime })\mapsto
\langle \sigma _{O}\mid \bar{\sigma}_{O}^{\prime }\rangle =(\sigma _{O}\mid
\sigma _{O}^{\prime })_{\mathcal{H}_{O}}\text{.} 
\]

In the following, we denote elements of $\mathcal{H}_{P\times O}$ by $\tau $
in order to simplify the notation. Each $\tau \in \mathcal{H}_{P\times O}=%
\mathcal{H}\otimes \overline{\mathcal{H}}_{O}$ gives rise to a continuous
linear map $\Theta _{\tau }:\mathcal{H}_{O}\rightarrow \mathcal{H}$ such
that, if $\tau =\sigma ^{1}\otimes \bar{\sigma}_{O}^{1}+...+\sigma
^{k}\otimes \bar{\sigma}_{O}^{k}$, for some $k\in \mathbb{N}$, where $\sigma
^{1},...,\sigma ^{k}$ and $\bar{\sigma}_{O}^{1},...,\bar{\sigma}_{O}^{k}$
are linearly independent sections of $L$ and $\bar{L}_{O}$, respectively,
then 
\begin{equation}
\Theta _{\tau }\sigma _{O}^{\prime }=\langle \sigma _{O}^{\prime }\mid \bar{%
\sigma}_{O}^{1}\rangle \sigma ^{1}+...+\langle \sigma _{O}^{\prime }\mid 
\bar{\sigma}_{O}^{k}\rangle \sigma ^{k}  \label{Theta}
\end{equation}%
for every $\sigma _{O}^{\prime }\in \mathcal{H}_{O}$.

Since $\mathcal{U}_{P\times O}=\mathcal{U}\otimes \overline{\mathcal{U}}_{O}$
it follows that, for every $g\in G$, $\ $%
\begin{eqnarray*}
&&(\mathcal{U}_{P\times O})_{g}(\sigma ^{1}\otimes \bar{\sigma}%
_{O}^{1}+...+\sigma ^{k}\otimes \bar{\sigma}_{O}^{k})= \\
&=&(\mathcal{U}_{g}\otimes (\overline{\mathcal{U}}_{O})_{g})(\sigma
^{1}\otimes \bar{\sigma}_{O}^{1}+...+\sigma ^{k}\otimes \bar{\sigma}_{O}^{k})
\\
&=&\mathcal{U}_{g}\sigma ^{1}\otimes (\overline{\mathcal{U}}_{O}^{G})_{g}%
\bar{\sigma}_{O}^{1}+...+\mathcal{U}_{g}\sigma ^{k}\otimes (\overline{%
\mathcal{U}}_{O}^{G})_{g}\bar{\sigma}_{O}^{k}\text{. }
\end{eqnarray*}%
Hence, $G$-invariance of $\tau =\sigma ^{1}\otimes \bar{\sigma}%
_{O}^{1}+...+\sigma ^{k}\otimes \bar{\sigma}_{O}^{k}$ means 
\[
\mathcal{U}_{g}\sigma ^{1}\otimes (\overline{\mathcal{U}}_{O})_{g}\bar{\sigma%
}_{O}^{1}+...+\mathcal{U}_{g}\sigma ^{k}\otimes (\overline{\mathcal{U}}%
_{O})_{g}\bar{\sigma}_{O}^{k}=\sigma _{P}^{1}\otimes \bar{\sigma}%
_{O}^{1}+...+\sigma _{P}^{k}\otimes \bar{\sigma}_{O}^{k} 
\]%
for all $g\in G$. We denote by $\mathcal{H}_{P\times O}^{G}$ the space of $G$%
-invariant elements of $\mathcal{H}_{P\times O}$. In other words, 
\[
\mathcal{H}_{P\times O}^{G}=\{\tau \in \mathcal{H}_{P\times O}\mid (\mathcal{%
U}_{P\times O})_{g}\tau =\tau \text{ for all }g\in G\}. 
\]

For $\tau \in \mathcal{H}_{P\times O}^{G},$ the map $\Theta _{\tau }$
intertwines the actions of $G$ on $\mathcal{H}_{O}$ and $\mathcal{H}.$ In
other words, for each $g\in G$,%
\begin{equation}
\Theta _{\tau }(\mathcal{U}_{O})_{g}=\mathcal{U}_{g}\Theta _{\tau }.
\label{intertwines}
\end{equation}%
If $\tau \neq 0$, then the restriction $\mathcal{U}_{\mid \mathrm{range~}%
\Theta _{\tau }}$ of the representation $\mathcal{U}$ to the range of $%
\Theta _{\tau }$ is a representation of $G$ on the range of $\Theta _{\tau }$
equivalent to $\mathcal{U}_{O}$. For $\tau ,\tau ^{\prime }\in \mathcal{H}%
_{P\times O}^{G}$, if the intersection of the range of $\Theta _{\tau }$ and
the range of $\Theta _{\tau ^{\prime }}$ is not empty$\neq 0$, then $\Theta
_{\tau }$ and $\Theta _{\tau ^{\prime }}$ have the same range and $\tau
^{\prime }$ is proportional to $\tau $. Thus, to each one-dimensional
subspace of $\mathcal{H}_{P\times O}^{G}$, generated by $\tau \neq 0$, there
corresponds an invariant subspace $\mathcal{H}_{\tau }=\mathrm{range}~\Theta
_{\tau }$ of $\mathcal{H}$ on which the quantization representation $%
\mathcal{U}$ is equivalent to $\mathcal{U}_{O}$.

Let $\Theta _{\tau }^{\ast }:\mathcal{H}\rightarrow \mathcal{H}_{O}$ denote
the adjoint of $\Theta _{\tau }$. For each $\sigma \in \mathcal{H}$ and $%
\sigma _{O}\in \mathcal{H}_{O}$, 
\[
(\Theta _{\tau }^{\ast }\sigma \mid \sigma _{O})_{\mathcal{H}_{O}}=(\sigma
\mid \Theta _{\tau }\sigma _{O})\text{,}
\]%
where $(\cdot \mid \cdot )$ is the scalar product in $\mathcal{H}$. Since
the quantization representations $\mathcal{U}_{O}$ $\ $and $\mathcal{U}$ are
unitary, equation (\ref{intertwines}) implies that%
\begin{eqnarray*}
((\mathcal{U}_{O})_{g}\Theta _{\tau }^{\ast }\sigma \mid \sigma _{O})_{%
\mathcal{H}_{O}} &=&(\Theta _{\tau }^{\ast }\sigma \mid (\mathcal{U}%
_{O})_{g^{-1}}\sigma _{O})_{\mathcal{H}_{O}}=(\sigma \mid \Theta _{\tau }(%
\mathcal{U}_{O})_{g^{-1}}\sigma _{O})= \\
&=&(\sigma \mid \mathcal{U}_{g^{-1}}\Theta _{\tau }\sigma _{O})=(\mathcal{U}%
_{g}\sigma \mid \Theta _{\tau }\sigma _{O})=(\Theta _{\tau }^{\ast }\mathcal{%
U}_{g}\sigma \mid \sigma _{O})_{\mathcal{H}_{O}}
\end{eqnarray*}%
for all $\sigma _{O}\in \mathcal{H}_{O}$, $\sigma \in \mathcal{H}$, and $%
g\in G$. Hence, $\Theta _{\tau }^{\ast }\mathcal{U}_{g}=(\mathcal{U}%
_{O})_{g}\Theta _{\tau }^{\ast }$. Therefore,%
\[
(\mathcal{U}_{O})_{g}\Theta _{\tau }^{\ast }\Theta _{\tau }=\Theta _{\tau
}^{\ast }\mathcal{U}_{g}\Theta _{\sigma }=\Theta _{\sigma }^{\ast }\Theta
_{\sigma }(\mathcal{U}_{O})_{g}.
\]%
Since $\mathcal{U}_{O}$ is irreducible, Schur's Lemma implies that the
mapping $\Theta _{\tau }^{\ast }\Theta _{\tau }$ $:\mathcal{H}%
_{O}\rightarrow \mathcal{H}_{O}$ is a multiple of identity on $\mathcal{H}%
_{O}$. Thus, $\Theta _{\tau }^{\ast }\Theta _{\tau }=\lambda _{\tau }\mathrm{%
Id}_{\mathcal{H}_{O}}$. On the other hand, 
\[
(\Theta _{\tau }\Theta _{\tau }^{\ast })^{2}=\Theta _{\tau }\Theta _{\tau
}^{\ast }\Theta _{\tau }\Theta _{\tau }^{\ast }=\Theta _{\tau }(\Theta
_{\tau }^{\ast }\Theta _{\tau })\Theta _{\tau }^{\ast }=\Theta _{\tau
}(\lambda _{\tau }\mathrm{Id}_{\mathcal{H}_{O}})\Theta _{\tau }^{\ast
}=\lambda _{\tau }\Theta _{\tau }\Theta _{\tau }^{\ast }\text{.}
\]%
If $\tau \neq 0$, then $\lambda _{\tau }>0$, and taking $\tau ^{\prime
}=(\lambda _{\tau })^{-1/2}\tau $, we get $\lambda _{\tau ^{\prime }}=1$,
and $(\Theta _{\tau ^{\prime }}\Theta _{\tau ^{\prime }}^{\ast })^{2}=\Theta
_{\tau ^{\prime }}\Theta _{\tau ^{\prime }}^{\ast }.$ For each $\sigma
,\sigma ^{\prime }\in \mathcal{H}$, we have 
\[
\left( \sigma \mid \Theta _{\tau ^{\prime }}\Theta _{\tau ^{\prime }}^{\ast
}\sigma ^{\prime }\right) =\left( \Theta _{\tau ^{\prime }}^{\ast }\sigma
\mid \Theta _{\tau ^{\prime }}^{\ast }\sigma ^{\prime }\right) =\left(
\Theta _{\tau ^{\prime }}\Theta _{\tau ^{\prime }}^{\ast }\sigma \mid \sigma
^{\prime }\right) ,
\]%
which implies that $\Theta _{\tau ^{\prime }}\Theta _{\tau ^{\prime }}^{\ast
}$ is self adjoint on $\mathcal{H}$. Thus, $\Theta _{\tau ^{\prime }}\Theta
_{\tau ^{\prime }}^{\ast }$ is the projection operator onto $\mathcal{H}%
_{\tau }=\mathcal{H}_{\tau ^{\prime }}$. In particular, it implies that $%
\mathcal{H}_{\tau ^{\prime }}$ is a closed subspace of $\mathcal{H}$. For
every $\sigma _{O}$ and $\sigma _{O}^{\prime }$ in $\mathcal{H}_{O}$ 
\[
\left( \sigma _{O}\mid \sigma _{O}^{\prime }\right) _{\mathcal{H}%
_{O}}=\left( \sigma _{O}\mid \mathrm{Id}_{\mathcal{H}_{O}}\sigma
_{O}^{\prime }\right) _{\mathcal{H}_{O}}=\left( \sigma _{O}\mid \Theta
_{\tau ^{\prime }}^{\ast }\Theta _{\tau ^{\prime }}\sigma _{O}^{\prime
}\right) _{\mathcal{H}_{O}}=\left( \Theta _{\tau ^{\prime }}\sigma _{O}\mid
\Theta _{\tau ^{\prime }}\sigma _{O}^{\prime }\right) \text{,}
\]%
which implies that $\Theta _{\tau ^{\prime }}:\mathcal{H}_{O}\rightarrow 
\mathcal{H}$ is an isometry. Thus, the restriction of $\mathcal{U}$ to $%
\mathcal{H}_{\tau ^{\prime }}$ is unitarily equivalent to $\mathcal{U}_{O}$.

Conversely, let $E$ be a closed $G$-invariant subspace of $\mathcal{H}$ such
that the restriction of $\mathcal{U}$ to $E$ is unitarily equivalent to $%
\mathcal{U}_{O}$. This means that there exists an isometry $\Theta :\mathcal{%
H}_{O}\rightarrow \mathcal{H}$ with $\mathrm{range}\Theta =E$ such that $%
\mathcal{U}\circ \Theta =\Theta \circ \mathcal{U}_{O}$. Since $\mathcal{H}%
_{P\times O}=\mathcal{H}\otimes \mathcal{\bar{H}}_{O}$ is naturally
isomorphic to $\mathrm{Hom}(\mathcal{H}_{O},\mathcal{H})$ and $\Theta $ is $%
G $-equivariant, it follows that there exists $\tau \in \mathcal{H}_{P\times
O}^{G}$ such that $\Theta =\Theta _{\tau }$. Hence, closed subspaces of $%
\mathcal{H}$ on which the restriction of $\mathcal{U}$ is unitarily
equivalent to $\mathcal{U}_{O}$ are completely determined by $\mathcal{H}%
_{P\times O}^{G}$.

Let $(\tau _{i})_{i\in \mathbb{N}}$ be a basis in $\mathcal{H}_{P\times
O}^{G}$. Without loss of generality, we may choose sections $\tau _{i}$ so
that operators $\Theta _{\tau _{i}}\Theta _{\tau _{i}}^{\ast }$ are
projections onto mutually orthogonal subspaces of $\mathcal{H}$. Then, 
\begin{equation}
\Pi _{O}=\dsum\limits_{i}\Theta _{\tau _{i}}\Theta _{\tau _{i}}^{\ast }
\label{PiO}
\end{equation}%
is a projection operator onto a subspace of $\mathcal{H}$ on which $\mathcal{%
U}$ is equivalent to a direct sum of $\mathcal{U}_{O}$ with multiplicity
given by the dimension of $\mathcal{H}_{P\times O}^{G}$ which may be
infinite.

By Theorem \ref{theorem3}, there is a natural isomorphism $\Xi $ of the
space $S_{F}^{\infty }(L_{P\times O})^{G}$ of $G$-invariant polarized
sections of $L_{P\times O}$ onto the representation space of quantization of
algebraic reduction at $O$ provided $J_{P\times O}^{-1}(0)$ contains a
Lagrangian manifold. Hence, the subspace $\Xi (\mathcal{H}_{P\times O}^{G})$
of the representation space of quantization of algebraic reduction at $O$
gives a complete description of the closed subspace of $\mathcal{H}$ on
which $\mathcal{U}$ is unitarily equivalent to a direct sum of $\mathcal{U}%
_{O}$ with multiplicity given by the dimension of $\mathcal{H}_{P\times
O}^{G}$.

We summarize our results below.

\begin{theorem}
\label{theorem5}Assume the following:

\noindent 1. $(P,\omega )$ is a quantizable K\"{a}hler manifold, $%
L\rightarrow P$ is a prequantization line bundle, and $F$ is the
distribution of antiholomorphic directions on $P.$ There is an action of a
connected Lie group $G$ on $P$ which preserves $F$ and has an $Ad^{\ast }$%
-equivariant momentum map $J:P\rightarrow \mathfrak{g}^{\ast }$. Moreover,
quantization of $(P,\omega )$ gives a unitary representation $\mathcal{U}$
of $G$ on the space $\mathcal{H}$ of square integrable holomorphic sections
of $L\rightarrow P$.

\noindent 2. $O$ is a quantizable coadjoint orbit of $G$ admitting a K\"{a}%
hler polarization $F_{O}$ such that quantization of $(O,\omega _{O})$ in
terms of the polarization $F_{O}$ gives rise to an irreducible unitary
representation $\mathcal{U}_{O}$ of $G$.

\noindent 3. There exists a Lagrangian subspace of $(P,\omega )$\ contained
in $J^{-1}(O)$.

Then, quantization of algebraic reduction at $O$, defines a projection
operator $\Pi _{O}$ on $\mathcal{H}$ with closed range on which the
quantization representation $\mathcal{U}$\ is unitarily equivalent to a
multiple of $\mathcal{U}_{O}$.
\end{theorem}

\end{document}